\documentclass[aos,MSNbibl,nameyear,seceqn,dvips]{arximspdf}
\usepackage{algorithm,algorithmic}
\usepackage{graphicx}

\doi{10.1214/12-AOS1012} 
\volume{40}
\issue{3}
\pubyear{2012}
\firstpage{1682}
\lastpage{1713}

\makeatletter

\newcommand{\rrvert}{\vert}
\newcommand{\llvert}{\vert}
\newcommand{\eqref}[1]{(\ref{#1})}
\renewcommand{\citep}[1]{\citeauthor{#1} \citeyear{#1}}
\newtheorem{reptheorem}{Theorem}
\newtheorem{theorem}{Theorem}
\newtheorem{lemma}{Lemma}
\newtheorem{proposition}{Proposition}

\newtheorem{problem}{Problem}

\newproclaim{definition}{Definition}
\newproclaim{example}{Example}
\newproclaim{remark}{Remark}

\newcommand{\dir}[1]{\stackrel{\rightarrow}{#1}}
\newcommand{\bidir}[1]{\stackrel{\leftrightarrow}{#1}}
\def\L{\Lambda}
\def\O{\Omega}
\def\S{\Sigma}
\def\det{\operatorname{det}}
\def\R{\mathbb{R}}
\def\AA{\mathbf{A}}
\def\bb{\mathbf{b}}
\def\JJ{\mathbf{J}}
\def\bidir{\leftrightarrow}
\def\dir{\rightarrow}
\def\eqbi{\mbox{$\,\circ$---$\circ\,$}}
\newcommand{\xleftarrow}{\leftarrow}
\newcommand{\xrightarrow}{\rightarrow}
\newcommand{\pa}{\mathrm{pa}} 
\newcommand{\sib}{\operatorname{sib}} 
\newcommand{\htr}{\operatorname{htr}} 

\newcommand{\bi}{\leftrightarrow}

\makeatother

\begin{document}
\begin{frontmatter}

\title{Half-trek criterion for generic identifiability of~linear structural equation models}
\runtitle{Identifiability of linear structural equation models}

\begin{aug}
\author[A]{\fnms{Rina} \snm{Foygel}\corref{}\ead[label=e1]{rina@galton.uchicago.edu}},
\author[B]{\fnms{Jan} \snm{Draisma}\ead[label=e2]{j.draisma@tue.nl}\thanksref{t1}}
\and
\author[A]{\fnms{Mathias} \snm{Drton}\ead[label=e3]{drton@uchicago.edu}\thanksref{t2}}

\thankstext{t1}{Supported by a Vidi grant from The Netherlands
Organisation for Scientific Research (NWO).}
\thankstext{t2}{Supported by the NSF under Grant DMS-07-46265 and by
an Alfred P. Sloan
Fellowship.}
\runauthor{R. Foygel, J. Draisma and M. Drton}
\affiliation{University of Chicago, Eindhoven University of Technology
and Centrum~voor~Wiskunde en Informatica, and University of Chicago}

\address[A]{R. Foygel\\
M. Drton\\
Department of Statistics\\
University of Chicago \\
Chicago, Illinois\\
USA\\
\printead{e1}\\
\printead{e3}}

\address[B]{J. Draisma\\
Department of Mathematics \\
\quad and Computer Science\\
Eindhoven University of Technology\\
Eindhoven\\
The Netherlands\\
and\\
Centrum voor Wiskunde en Informatica\\
Amsterdam\\
The Netherlands\\
\printead{e2}}


\end{aug}

\received{\smonth{7} \syear{2011}}
\revised{\smonth{2} \syear{2012}}

%
\begin{abstract}
A linear structural equation model relates random variables of
interest and corresponding Gaussian noise terms via a linear
equation system. Each such model can be represented by a mixed
graph in which directed edges encode the linear equations and
bidirected edges indicate possible correlations among noise terms.
We study parameter identifiability in these models, that is, we ask
for conditions that ensure that the edge coefficients and
correlations appearing in a~linear structural equation model can be
uniquely recovered from the covariance matrix of the associated
distribution. We treat the case of generic identifiability, where
unique recovery is possible for almost every choice of parameters.
We give a new graphical condition that is sufficient for generic
identifiability and can be verified in time that is polynomial in
the size of the graph. It improves criteria from prior work and
does not require the directed part of the graph to be acyclic. We
also develop a related necessary condition and examine the ``gap''
between sufficient and necessary conditions through simulations on graphs
with $25$ or $50$ nodes,
as well as exhaustive algebraic computations for graphs with up to five
nodes.
\end{abstract}

%


\end{frontmatter}

\section{Introduction}
\label{secintroduction}

When modeling the joint distribution of a random vector
$X=(X_1,\ldots,X_m)^T$, it is often natural to appeal to noisy
functional relationships. In other words, each variable $X_w$ is
assumed to be a function of the remaining variables and a stochastic
noise term $\varepsilon_w$. The resulting models are known as linear
structural equation
models when the relationship is linear, that is, when
%
\begin{equation}
\label{eqsem} X_w = \lambda_{0w}+ \sum
_{v\not=w} \lambda_{vw} X_v+
\varepsilon_w,\qquad w=1,\ldots,m,
\end{equation}
or in vectorized form with a matrix $\Lambda=(\lambda_{vw})$ that is
tacitly assumed to have zeros along the diagonal,
%
\begin{equation}
\label{eqsem-vec} X = \lambda_0 + \Lambda^T X +
\varepsilon.
\end{equation}
The classical distributional assumption is that the error
vector $\varepsilon=(\varepsilon_1,\ldots,\break\varepsilon_m)^T$ has a
multivariate normal distribution with zero mean and some covariance
matrix $\Omega=(\omega_{vw})$. Writing $I$ for the identity matrix,
it follows that~$X$ has a~multivariate normal distribution with mean
vector $(I-\Lambda)^{-T}\lambda_0$ and covariance matrix
%
\begin{equation}
\label{eqcov-mx} \Sigma= (I-\Lambda)^{-T}\Omega(I-\Lambda)^{-1}.
\end{equation}
Background on structural equation modeling can be found, for instance,
in \citet{bollen1989}. As emphasized in \citet{spirtes2000} and \citet
{pearl2000}, their
great popularity in applied sciences is due to the natural causal
interpretation of the involved functional relationships.

Interesting models are obtained by imposing some pattern of zeros
among the coefficients $\lambda_{vw}$ and the covariances
$\omega_{vw}$. It is convenient to think of the zero patterns as
being associated with a mixed graph that contains directed edges $v\to
w$ to indicate possibly nonzero coefficients $\lambda_{vw}$, and
bidirected edges $v\bi w$ when $\omega_{vw}$ is a possibly nonzero
covariance; in figures we draw the bidirected edges dashed for better
distinction. Mixed graph representations have first been advocated in
\citeauthor{wright1921} (\citeyear{wright1921,wright1934}) and are also
known as path diagrams.
We briefly illustrate this in the next example, which gives the simplest
version of what are often referred to as instrumental variable models;
see also \mbox{\citet{didelez2010}}.

\begin{figure}

\includegraphics{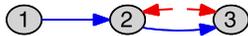}

\caption{Mixed graph for the instrumental variable model.}
\label{figiv-graph}
\end{figure}

%
\begin{example}[(IV)]
\label{exiv}
Suppose that, as in \citet{evans1999}, we record an infant's birth
weight ($X_3$), the level of maternal smoking during pregnancy
($X_2$) and the cigarette tax rate that applies ($X_1$). A model
of interest, with mixed graph in Figure~\ref{figiv-graph}, assumes
\[
X_1 = \lambda_{01}+\varepsilon_1,\qquad
X_2 = \lambda_{02}+\lambda_{12}X_1
+ \varepsilon_2,\qquad
X_3 = \lambda_{03}+\lambda_{23} X_2
+ \varepsilon_3,
\]
with an error vector $\varepsilon$ that has zero mean
vector and covariance matrix
\begin{eqnarray*}
\Omega= \pmatrix{ \omega_{11} & 0&0\vspace*{2pt}
\cr
0&
\omega_{22} & \omega_{23}\vspace*{2pt}
\cr
0&
\omega_{23} & \omega_{33}}.
\end{eqnarray*}
The possibly nonzero entry $\omega_{23}$ can absorb the effects
that unobserved confounders (such as age, income, genetics, etc.)
may have on both $X_2$ and~$X_3$; compare
\citet{richardson2002} and \citet{wermuth2010} for background on
mixed graph
representations of latent variable problems.
\end{example}

Formally, a mixed graph is a triple $G=(V,D,B)$, where $V$ is a finite
set of nodes and $D,B\subseteq V\times V$ are two sets of edges. In
our context, the nodes correspond to the random variables
$X_1,\ldots,X_m$, and we simply let $V=[m]:=\{1,\ldots,m\}$. The pairs
$(v,w)$ in the set $D$ represent directed edges and we will always
write $v\to w$; $v\to w\in D$ does not
imply $w\to v\in D$. The pairs in $B$ are bidirected edges $v\bi w$;
they have no orientation, that is, $v\bi w\in B$ if and only if $w\bi
v\in B$.
Neither the bidirected part $(V,B)$ nor the directed part $(V,D)$
contain self-loops, that is, $v\to v\notin D$ and $v\bi v\notin B$
for all $v\in V$. If the directed part $(V,D)$ does not contain
directed cycles (i.e., no cycle $v\to\cdots\to v$ can be formed from
the edges in $D$), then the mixed graph $G$ is said to be acyclic.

Let $\mathbb{R}^D$ be the set of real $m\times m$-matrices
$\Lambda=(\lambda_{vw})$ with support $D$, that is, $\lambda_{vw}=0$
if $v\to w\notin D$. Write $\mathbb{R}^D_\mathrm{reg}$ for the
subset of matrices $\Lambda\in\mathbb{R}^D$ for which $I-\Lambda$ is
invertible, where $I$ denotes the identity matrix. [If $G$ is
acyclic, then $\mathbb{R}^D=\mathbb{R}^D_\mathrm{reg}$; see the remark
after equation~\eqref{eqtrek-rule}.] Similarly,
let $\operatorname{PD}_m$ be the cone of positive definite symmetric
$m\times m$-matrices $\Omega=(\omega_{vw})$ and define
$\operatorname{PD}(B)\subset\operatorname{PD}_m$ to be the subcone of matrices
with support $B$, that is, $\omega_{vw}=0$ if $v\not= w$ and $v\bi
w\notin B$.

%
\begin{definition}
\label{defsem}
The linear structural equation model given by a mixed graph
$G=(V,D,B)$ on $V=[m]$ is the family of all $m$-variate normal
distributions with covariance
matrix
\[
\Sigma= (I-\Lambda)^{-T}\Omega(I-\Lambda)^{-1}
\]
for $\Lambda\in\mathbb{R}^D_\mathrm{reg}$ and
$\Omega\in\operatorname{PD}(B)$.
\end{definition}

The first question that arises when specifying a linear structural
equation model is whether the model is identifiable in the sense that
the parameter matrices $\Lambda\in\mathbb{R}^D_\mathrm{reg}$ and
$\Omega\in\operatorname{PD}(B)$ can be uniquely recovered from the normal
distribution they define. Clearly, this is equivalent to asking
whether they can be recovered from the distribution's covariance
matrix, and thus we ask whether the
\textit{fiber}
%
\begin{equation}
\label{eqfiber} \mathcal{F}(\Lambda,\Omega) = \bigl\{ \bigl(
\Lambda',\Omega'\bigr)\in\Theta\dvtx
\phi_G\bigl(\Lambda',\Omega'\bigr)=
\phi_G(\Lambda,\Omega) \bigr\}
\end{equation}
is equal to $\{(\Lambda,\Omega)\}$. Here, we introduced the shorthand
$\Theta:= \mathbb{R}^D_\mathrm{reg}\times\operatorname{PD}(B)$. Put
differently, identifiability holds if the parametrization map
%
\begin{equation}
\label{eqphiG} \phi_G\dvtx(\Lambda,\Omega)\mapsto(I-
\Lambda)^{-T}\Omega(I-\Lambda)^{-1}
\end{equation}
is injective on $\Theta$, or a suitably large subset.

%
\begin{example}[(IV, continued)]
\label{exivcont}
In the instrumental variable model associated with the graph in
Figure~\ref{figiv-graph},
\begin{eqnarray*}
\Sigma&=& (\sigma_{vw}) \\
&=& \pmatrix{ 1 & -\lambda_{12} & 0
\vspace*{2pt}
\cr
0 & 1 & -\lambda_{23}\vspace*{2pt}
\cr
0 & 0 &
1}^{-T} \pmatrix{ \omega_{11} & 0&0\vspace*{2pt}
\cr
0&
\omega_{22} & \omega_{23}\vspace*{2pt}
\cr
0&
\omega_{23} & \omega_{33}} \pmatrix{ 1 & -
\lambda_{12} & 0\vspace*{2pt}
\cr
0 & 1 & -\lambda_{23}
\vspace*{2pt}
\cr
0 & 0 & 1}^{-1}
\\
& =& \pmatrix{ \omega_{11} & \omega_{11}\lambda_{12}
& \omega_{11}\lambda_{12} \lambda_{23}
\vspace*{2pt}
\cr
\omega_{11}\lambda_{12} &
\omega_{22} + \omega_{11}\lambda_{12}^2
& \omega_{23} + \lambda_{23}\sigma_{22}
\vspace*{2pt}
\cr
\omega_{11}\lambda_{12}\lambda_{23}
& \omega_{23} + \lambda_{23}\sigma_{22}&
\omega_{33} + 2\omega_{23}\lambda_{23} +
\lambda_{23}^2\sigma_{22}}.
\end{eqnarray*}
Despite the presence of both the edges $2\to3$ and $2\bi3$, we can
recover $\Lambda$ (and thus also $\Omega$) from $\Sigma$ using that
\[
\lambda_{12} = \frac{\sigma_{12}}{\sigma_{11}}, \qquad\lambda_{23} =
\frac{\sigma_{13}}{\sigma_{12}}.
\]
The first denominator $\sigma_{11}$ is always positive since
$\Sigma$ is positive definite. The second denominator $\sigma_{12}$
is zero if and only if $\lambda_{12}=0$. In other words, if the
cigarette tax ($X_1$) has no effect on maternal smoking during
pregnancy ($X_2$), then there is no way to distinguish between the
causal effect of smoking on birth weight (coefficient
$\lambda_{23}$) and the effects of confounding variables
(error covariance $\omega_{23}$). Indeed, the map $\phi_G$ is injective
only on the subset of $\Theta$ with $\lambda_{12}\not=0$.
\end{example}

In this paper we study the kind of identifiability encountered in the
instrumental variables example. The statistical literature often
refers to this as almost-everywhere identifiability to express that
the exceptional pairs $(\Lambda,\Omega)$ with fiber cardinality
$|\mathcal{F}(\Lambda,\Omega)|>1$ form a set of measure zero.
However, since the map~$\phi_G$ is rational, the exceptional sets are
well-behaved null sets, namely, they are algebraic subsets. An
algebraic subset $V\subset\Theta$ is a~subset that can be defined by
polynomial equations, and it is a proper subset of the open set
$\Theta$ unless it is defined by the zero polynomial. A~proper
algebraic subset has smaller dimension than $\Theta$ [see
\citet{cox2007}], and thus also measure zero; statistical work often
quotes the lemma in \citet{okamoto1973} for the latter fact. These
observations motivate the following definition and problem.

%
\begin{definition}
The mixed graph $G$ is said to be generically identifiable if
$\phi_G$ is injective on the complement $\Theta\setminus V$ of a
proper (i.e., strict) algebraic subset $V\subset\Theta$.
\end{definition}

%
\begin{problem}
\label{probcharacterize}
Characterize the mixed graphs $G$ that are generically identifiable.
\end{problem}

Despite the long history of linear structural equation models, the
problem just stated remains open, even when restricting to acyclic
mixed graphs. However, in the last two decades a number of graphical
conditions have been developed that are sufficient for generic
identifiability. We refer the reader, in particular, to
\citet{pearl2000}, \citeauthor{brito2002} (\citeyear{brito2002,brito2006}),
\citet{tian2009} and \citet{chan2010}, which each contain many
further references. To our knowledge, the condition that is of the most
general nature and most in the spirit of attempting to solve
Problem~\ref{probcharacterize} is the G-criterion of
\citet{brito2006}. This criterion, and in fact all other mentioned
work, uses linear algebraic techniques to solve the parametrized
equation systems that define the fibers $\mathcal{F}(\Lambda,\Omega)$.
Therefore, the G-criterion is in fact sufficient for the following stronger
notion of identifiability, which we have seen to hold for the graph
from Figure~\ref{figiv-graph}; recall the formulas given in
Example~\ref{exivcont}.

%
\begin{definition}
\label{defrational-ident}
The mixed graph $G$ is said to be rationally identifiable if there
exists a proper algebraic subset $V\subset\Theta$ and a rational map
$\psi$ such that $\psi\circ\phi_G(\Lambda,\Omega) =
(\Lambda,\Omega)$ for all $(\Lambda,\Omega)\in\Theta\setminus V$.
\end{definition}

The main results of our paper give a graphical condition that is
sufficient for rational identifiability and that is strictly stronger
than the G-criterion of \citet{brito2006} when applied to acyclic
mixed graphs. Moreover, the new condition, which we name the half-trek
criterion, also applies to cyclic graphs, for which little prior
work exists. The approach we take also yields a necessary condition, or,
more precisely put, a graphical condition that is sufficient for $G$ (or
rather the map $\phi_G$) to be \textit{generically infinite-to-one}. That
is, the condition implies that the fiber $\mathcal{F}(\Lambda,\Omega)$
is infinite for all pairs $(\Lambda,\Omega)$ outside a proper algebraic
subset of~$\Theta$. Hardly any previous work on such ``negative''
graphical conditions seems to exist.
Our main results just described are stated in detail in
Section~\ref{secmain-results} and proven in Section~\ref{secproofs-H}
and in Sections~2 and~3 of the Supplementary Material [\citet{Supplement}].
The comparison to the G-criterion
is made in Section~\ref{secG-crit}, with some proofs deferred to
Section~4
of the supplement. Some interesting examples are visited
in Section~\ref{secexamples}. Those include examples that do not seem
to be covered by any known graphical criterion.

A major motivation for this paper is the \textit{complexity} of deciding
whether a~given graph is rationally identifiable. In \citet{garcia2010}
this question is proved to be decidable using computational algebraic
geometry, and in Section~8 of the supplement we
give a variant of that
approach in which the size of the input to Buchberger's algorithm is
significantly reduced. However, there is no reason to believe that
this approach yields an algorithm whose running time is bounded by some
polynomial in the size of the input, namely, the mixed graph $G$. Faced
with this situation, one naturally wonders whether this decision problem
is at all contained in complexity class NP, which requires that for all
rationally identifiable $G$ there exists a certificate for rational
identifiability that can be checked in polynomial time. This is by
no means clear to us. For instance, while in Example~\ref{exivcont}
the rational inverse map of the parametrization happens to be rather
small in terms of bit-size, it is unclear why for general rationally
identifiable~$G$ there should be a rational map that can in polynomial
time be checked to be inverse to the parametrization (on the other hand,
there is no reason why efficiently checkable certificates would have
to be of this form). By contrast,
our half-trek criteria for rational identifiability and for being
generically infinite-to-one turn out not only to have efficiently
checkable certificates for positive instances (which will be evident
from the criteria's definitions) but even to be in complexity class
$\mathrm{P}\subseteq\mathrm{NP}$. Indeed, in Section~\ref{secalgorithms}
we develop polynomial-time algorithms for checking our graphical
conditions from Section~\ref{secmain-results}, and correctness of
those algorithms is proven in Section~6 of the
supplement.

The examples shown in Section~\ref{secexamples} were found as part of
an exhaustive study of the identifiability properties of all mixed
graphs with up to 5 nodes, in which we compare the aforementioned,
generally applicable but inefficient techniques from computational
algebraic geometry with our half-trek criteria. The results of these
computations are given in Section~\ref{seccomputations}. That section
further contains, as proof of concept, the result of simulations for
graphs on 25 or 50 nodes, based on the polynomial-time algorithms from
Section~\ref{secalgorithms}. Finally, in
Section~\ref{secgraph-decompositions}, we describe how our half-trek
methods behave with respect to a graph decomposition technique for
acyclic mixed graphs that is due to \citet{tian2005}; somewhat
surprisingly, this leads to a strengthening of our sufficient
condition. Concluding remarks are given in
Section~\ref{secconclusion}.


\section{Preliminaries on treks}

A path from node $v$ to node $w$ in a mixed graph $G=(V,D,B)$ is a
sequence of edges, each from either $D$ or $B$, that connect the
consecutive nodes in a sequence of nodes beginning at $v$ and ending
in $w$. We do not require paths to be simple or even to obey
directions, that is, a path may
include a particular edge more than once, the nodes that are part
of the edges need not all be distinct, and directed edges may be
traversed in the wrong direction. A path $\pi$ from $v$ to $w$
is a \textit{directed path} if all its edges are directed and pointing
to $w$, that is, $\pi$ is of the form
\[
v=v_0\to v_1\to\cdots\to v_r=w.
\]

In a covariance matrix in a structural equation model, that is, a
matrix structured as in Definition~\ref{defsem}, the entry
$\sigma_{vw}$ is a sum of terms that correspond to certain paths from
$v$ to $w$. For instance, in Example~\ref{exivcont}, the variance
%
\begin{equation}
\label{eqex-trekrule-sigma33} \sigma_{33} = \omega_{33} +
\omega_{23}\lambda_{23} +\omega_{23}
\lambda_{23} + \lambda_{23}^2\omega_{22}
+ \lambda_{23}^2\lambda_{12}^2
\omega_{11}
\end{equation}
is a sum of five terms that are associated, respectively, with the
trivial path~$3$,
which has no edges,
and the four additional paths
\[
3\leftrightarrow2\rightarrow3,\qquad 3\leftarrow2\leftrightarrow3, \qquad3\leftarrow2
\rightarrow3,\qquad 3\leftarrow2\leftarrow1\rightarrow2\rightarrow3.
\]
In the literature, the paths that contribute to a covariance are known
as \textit{treks}; compare, for example, \citet{sullivant2010} and the
references
therein. A \textit{trek} from \textit{source} $v$ to \textit{target} $w$
is a path
from $v$ to $w$ whose consecutive edges do not have any colliding
arrowheads. In other words, a trek from $v$ to $w$ is a path of one of
the two
following forms:
\[
v^{\mathrm{L}}_{l}\leftarrow v^{\mathrm{L}}_{l-1}
\leftarrow\cdots\leftarrow v^{\mathrm{L}}_1\leftarrow
v^{\mathrm{L}}_0\longleftrightarrow v^{\mathrm{R}}_0
\dir v^{\mathrm{R}}_1\dir\cdots\dir v^{\mathrm{R}}_{r-1}
\dir v^{\mathrm{R}}_{r}
\]
or
\[
v^{\mathrm{L}}_{l}\leftarrow v^{\mathrm{L}}_{l-1}
\leftarrow\cdots\leftarrow v^{\mathrm{L}}_1\xleftarrow v^{\mathrm
{T}}\xrightarrow
v^{\mathrm{R}}_1\dir\cdots\dir v^{\mathrm{R}}_{r-1}\dir
v^{\mathrm{R}}_{r},
\]
where the endpoints are $v^{\mathrm{L}}_{l}=v$, $v^{\mathrm{R}}_{r}=w$.
In the first case, we say that the left-hand side of $\pi$, written
$\operatorname{Left}({\pi})$, is the set of nodes
$\{v^{\mathrm{L}}_0,v^{\mathrm{L}}_1,\ldots,v^{\mathrm{L}}_{l}\}$,
and the
right-hand side, written $\operatorname{Right}({\pi})$, is the set of nodes
$\{v^{\mathrm{R}}_0,v^{\mathrm{R}}_1,\ldots,v^{\mathrm{R}}_{r}\}$.
In the
second case,
$\operatorname{Left}({\pi})=\{v^{\mathrm{T}},v^{\mathrm{L}}_1,\ldots
,v^{\mathrm
{L}}_{l}\}$,
and $\operatorname{Right}({\pi})=\{v^{\mathrm{T}},v^{\mathrm{R}}_1,\ldots
,v^{\mathrm
{R}}_{r}\}$---note that the \textit{top} node $v^{\mathrm{T}}$ is part of both
sides of
the trek. As pointed out before, paths and, in particular, treks are
not required to be simple. A~trek $\pi$ may thus pass through a node
on both its left- and right-hand sides. If the graph contains a
cycle, then the left- or right-hand side of $\pi$ may contain this
cycle. Any directed path is a trek; in this case
$|\operatorname{Left}({\pi})|=1$ or $|\operatorname{Right}({\pi})|=1$ depending
on the direction in
which the path is traversed. A~trek from $v$ to
$v$ may have no edges, in which case $v$ is the top node, and
$\operatorname{Left}({\pi})=\operatorname{Right}({\pi})=\{v\}$, and we call the
trek trivial.

A trek is therefore obtained by concatenating two directed paths at a
common top
node or by joining them with a bidirected edge, and the connection
between the matrix entries and treks is due to the fact that
%
\begin{equation}
\label{eqdirected-paths} \bigl((I- \Lambda)^{-1}\bigr)_{vw}
= \sum_{\pi\in\mathcal{P}(v,w)} \prod_{x\rightarrow y \in\pi}
\lambda_{xy},
\end{equation}
where $\mathcal{P}(v,w)$ is the set of directed paths from $v$ to $w$
in $G$. The equality in~(\ref{eqdirected-paths}) follows by writing
$(I-\L)^{-1}=I+\L+\L^2+\cdots$. For a precise statement about the form
of the covariance matrix $\Sigma$, let $\mathcal{T}(v,w)$ be the set
of all treks from $v$ to~$w$. For a trek $\pi$ that contains no
bidirected edge and has top node $v$, define a trek monomial as
\[
\pi(\lambda,\omega) = \omega_{vv}\prod_{x\to y\in\pi}
\lambda_{xy}.
\]
For a trek $\pi$ that contains a bidirected edge $v\bi w$, define
the trek monomial as
\[
\pi(\lambda,\omega) = \omega_{vw}\prod_{x\to y\in\pi}
\lambda_{xy}.
\]
The following rule [\citet{spirtes2000}, \citeauthor{wright1921}
(\citeyear{wright1921,wright1934})]
expresses the covariance matrix $\Sigma$ as a summation over treks;
compare the example in~(\ref{eqex-trekrule-sigma33}).
%
\begin{trekrule*} The covariance matrix $\Sigma$ for a mixed graph $G$
is given by
%
\begin{equation}
\label{eqtrek-rule} \sigma_{vw} = \sum_{\pi\in\mathcal{T}(v,w)}
\pi(\lambda,\omega);
\end{equation}
\end{trekrule*}

If $G$ is acyclic, then $\L^k=0$ for all $k\ge m$, and so the
expression in~(\ref{eqdirected-paths}) is polynomial. Similarly,
(\ref{eqtrek-rule}) writes $\sigma_{vw}$ as a polynomial.
If $G$ is cyclic, then one obtains power series that converge if the
entries of $\L$ are small enough. However, in the proofs of
Section~\ref{secproofs-H} it will also be useful to treat these as
formal power series.

Our identifiability results involve conditions that refer to paths
that we term half-treks. A \textit{half-trek} $\pi$ is a trek with
$\llvert\operatorname{Left}({\pi})\rrvert=1$, meaning that $\pi$ is of the
form
\[
v^{\mathrm{L}}_0\bidir{} v^{\mathrm{R}}_0\dir
v^{\mathrm{R}}_1\dir\cdots\dir v^{\mathrm{R}}_{r-1}\dir
v^{\mathrm{R}}_{r}
\]
or
\[
v^{\mathrm{T}}\dir v^{\mathrm{R}}_1\dir\cdots\dir
v^{\mathrm{R}}_{r-1}\dir v^{\mathrm{R}}_{r}.
\]

\begin{figure}[b]

\includegraphics{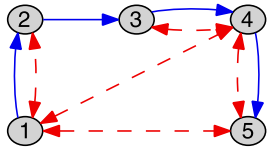}

\caption{An acyclic mixed graph.}\label{figtrekexamples}
\end{figure}

%
\begin{example} In the graph shown in Figure~\ref{figtrekexamples},
\begin{longlist}[(a)]
\item[(a)] neither $\pi_1\dvtx2\dir3\dir4\leftarrow3$ nor $\pi_2\dvtx
3\dir4\bidir1$
are treks, due to the colliding arrowheads at node $4$.
\item[(b)] $\pi\dvtx2\leftarrow1 \bidir4 \dir5$ is a trek, but not a
half-trek. $\operatorname{Left}({\pi})=\{1,2\}$ and $\operatorname{Right}({\pi
})=\{4,5\}$.
\item[(c)] $\pi\dvtx1\dir2\dir3$ is a half-trek with $\operatorname
{Left}({\pi})=\{1\}
$ and $\operatorname{Right}({\pi})=\{1,2,3\}$.
\end{longlist}
\end{example}

It will also be important to consider sets of treks. For a set of $n$
treks, $\Pi=\{\pi_1,\ldots,\pi_n\}$, let $x_i$ and $y_i$ be the source
and the target of $\pi_i$, respectively. If the sources are all
distinct, and the targets are all distinct, then we say that $\Pi$ is
a \textit{system of treks} from $X=\{x_1,\ldots,x_n\}$ to
$Y=\{y_1,\ldots,y_n\}$, which we write as $\Pi\dvtx X\rightrightarrows Y$.
Note that there may be overlap
between the sources in $X$ and the targets in $Y$, that is, we might
have $X\cap
Y\not=\varnothing$. The system~$\Pi$ is a system\vadjust{\goodbreak} of half-treks if every
trek~$\pi_i$ is a half-trek. Finally, a set of treks
$\Pi=\{\pi_1,\ldots,\pi_n\}$ has no sided intersection if
\[
\operatorname{Left}({\pi_i})\cap\operatorname{Left}({\pi_j})=
\varnothing=\operatorname{Right}({\pi_i})\cap\operatorname{Right}({
\pi_j}) \qquad\forall i\neq j .
\]

%
\begin{example} Consider again the graph from Figure~\ref{figtrekexamples}.
\begin{longlist}[(a)]
\item[(a)] The pair of treks
\[
\pi_1\dvtx3\dir4\dir5, \qquad \pi_2\dvtx4\bidir1
\]
forms a system of treks $\Pi=\{\pi_1,\pi_2\}$ between $X=\{3,4\}$ and
$Y=\{1,5\}$. The node $4$ appears in both treks, but is in only the
right-hand side of $\pi_1$ and only the left-hand side of $\pi_2$.
Therefore, $\Pi$ has no sided intersection.
\item[(b)] The set $\Pi=\{\pi_1,\pi_2\}$ comprising the two treks
\[
\pi_1\dvtx1\bidir4, \qquad \pi_2\dvtx3\dir4\dir5
\]
is a system of treks between $X=\{1,3\}$ and $Y=\{4,5\}$. Since
node $4$ is in $\operatorname{Right}({\pi_1})\cap\operatorname{Right}({\pi_2})$,
the system $\Pi$ has a sided intersection.
\end{longlist}
\end{example}

\section{Main identifiability and nonidentifiability results}
\label{secmain-results}

Define the set of \textit{parents} of a node $v\in V$ as $\pa(v) = \{
w\dvtx
w\to v\in D\}$ and the set of \textit{siblings} as $ \sib(v)=\{w\dvtx
w\leftrightarrow v\in B\}$. Let $ \htr(v)$ be the set of nodes in
$V\setminus(\{v\}\cup\sib(v))$ that can be reached from $v$ via a
half-trek. These half-treks contain at least one directed edge. Put
differently, a node $w\not=v$ that is not a sibling of $v$ is in $
\htr(v)$ if $w$ is a proper descendant of $v$ or one of its siblings.
Here, the term \textit{descendant} refers to a node that can be reached
by a directed path.

%
\begin{definition}
\label{defhalf-trek-crit}
A set of nodes $Y\subset V$ satisfies the half-trek criterion with
respect to node $v\in V$ if
\begin{longlist}[(iii)]
\item[(i)] $|Y|=|\pa(v)|$,
\item[(ii)] $Y\cap(\{v\}\cup\sib(v) ) = \varnothing$, and
\item[(iii)] there is a system of half-treks with no sided
intersection from $Y$ to $\pa(v)$.
\end{longlist}
\end{definition}

We remark that if $\pa(v)=\varnothing$, then $Y=\varnothing$ satisfies
the half-trek criterion with respect to $v$. We are now ready to
state the main results of this paper.

%
\begin{theorem}[(HTC-identifiability)]\label{mainthm}
Let $(Y_v\dvtx v\in V)$ be a family of subsets of the vertex
set $V$ of a mixed graph $G$. If, for each node $v$, the set
$Y_v$ satisfies the half-trek criterion with respect to $v$, and
there is a total ordering $\prec$ on the vertex set $V$ such that
$w\prec v$ whenever $w\in Y_v\cap\htr(v)$, then $G$ is
rationally identifiable.
\end{theorem}

The existence of such a total ordering is equivalent to the relation
$w \in Y_v \cap\htr(v)$ not admitting cycles; given the family
$(Y_v\dvtx v \in V)$, this can clearly be tested in polynomial time in the
size of the graph. More importantly, as we show in
Section~\ref{secalgorithms}, HTC-identifiability itself can be
checked in polynomial\vadjust{\goodbreak} time. In that section we will also show that
the same is true for the following nonidentifiability criterion.

%
\begin{theorem}[(HTC-nonidentifiability)]\label{mainthm2}
Suppose $G$ is a mixed graph in which every family $(Y_v\dvtx v\in
V)$ of subsets of the vertex set $V$ either contains a~set $Y_v$
that fails to satisfy the half-trek criterion with respect to $v$ or
contains a pair of sets $(Y_v,Y_w)$ with $v\in Y_w$ and
$w\in Y_v$. Then the parametrization $\phi_G$ is generically
infinite-to-one.
\end{theorem}

The main ideas underlying the two results are as follows. Under the
conditions given in Theorem~\ref{mainthm}, it is possible to recover
the entries in the matrix~$\Lambda$, column-by-column, following the
given ordering of the nodes. Each column is found by solving a linear
equation system that can be proven to have a~unique solution. The
details of these computations are given in Section~\ref{secproofs-H},
where we prove Theorem~\ref{mainthm}. The proof of
Theorem~\ref{mainthm2} is also in Section~\ref{secproofs-H} and rests
on the fact that under the given conditions the Jacobian of $\phi_G$
cannot have full rank.

In light of the two theorems, we refer to a mixed graph $G$ as follows:
\begin{longlist}[(iii)]
\item[(i)] HTC-identifiable, if it satisfies the conditions of
Theorem~\ref{mainthm},
\item[(ii)] HTC-infinite-to-one, if it satisfies the conditions of
Theorem~\ref{mainthm2},
\item[(iii)] HTC-classifiable, if it is either HTC-identifiable or
HTC-infinite-to-one,
\item[(iv)] HTC-inconclusive, if it is not HTC-classifiable.
\end{longlist}
We now give a first example of an HTC-identifiable graph. Additional
examples will be given in Section~\ref{secexamples}, where we will
see graphs that are generically $h$-to-one with $2\le h<\infty$, but
also that HTC-inconclusive graphs may be rationally identifiable or
generically infinite-to-one.

%
\begin{example}
\label{exGvshalfexample-HTCpart}
The graph in Figure~\ref{figtrekexamples} is
HTC-identifiable, which can be shown as follows. Let
\[
Y_1=\varnothing,\qquad Y_2=\{5\},\qquad Y_3=\{2\},\qquad
Y_4=\{2\}, \qquad Y_5=\{3\} .
\]
Then each $Y_v$ satisfies the half-trek criterion with respect
to $v$ because,
\begin{longlist}[(a)]
\item[(a)] trivially, $\pa(v)=\varnothing$ for $v=1$;
\item[(b)] for $v=2$, we have $5\bidir1\dir2$;
\item[(c)] for $v=3$, we have $2\dir3$;
\item[(d)] for $v=4$, we have $2\dir3\dir4$; and
\item[(e)] for $v=5$, we have $3\dir4\dir5$.
\end{longlist}
Considering the descendant sets $ \htr(v)$, we find that
\begin{eqnarray*}
Y_1\cap\htr{(1)}&=&\varnothing,\qquad Y_2\cap\htr{(2)}=\{5\} ,\qquad
Y_3\cap\htr{(3)}=\varnothing,
\\
Y_4\cap\htr{(4)}&=&\{2\} , \qquad Y_5\cap\htr{(5)}=\{3\} .
\end{eqnarray*}
Hence, any ordering $\prec$ respecting $3\prec5\prec2\prec4$ will
satisfy the
conditions of Theorem~\ref{mainthm}.
\end{example}

A mixed graph $G=(V,D,B)$ is simple if there is at most one edge
between any pair of nodes, that is, if $D\cap B=\varnothing$ and $v\to w
\in D$ implies $w\to v\notin D$. As observed in
\citet{brito2002}, simple acyclic mixed graphs are rationally
identifiable; compare also Corollary 3 in \citet{drton2011}.
It is not difficult to see that Theorem~\ref{mainthm} includes this
observation as a special case.

%
\begin{proposition}
\label{propsimple-graphs}
If $G$ is a simple acyclic mixed graph, then $G$ is HTC-identifiable.
\end{proposition}
\begin{pf}
Since $G$ is simple, it holds for every node $v\in V$ that
$\pa(v)\cap\sib(v)=\varnothing$ and, thus, $\pa(v)$ satisfies the
half-trek criterion with respect to $v$. An acyclic graph has at
least one topological ordering $\prec$, that is, an ordering such
that $v\to w\in D$ only if $v\prec w$. In other words, $w\in\pa(v)$ implies
$w\prec v$. Hence, the family $(\pa(v)\dvtx v\in V)$ together with
a topological ordering $\prec$ satisfies the conditions of
Theorem~\ref{mainthm}.
\end{pf}

Another straightforward observation is that the map $\phi_G$ cannot be
generically finite-to-one if the dimension of the domain of definition
$\mathbb{R}^D_\mathrm{reg}\times\operatorname{PD}(B)$ is larger than the
space of $m\times m$ symmetric matrices that contains the image of
$\phi_G$. This occurs if $|D|+|B|$ is larger than ${{m}\choose{2}}$.
Theorem~\ref{mainthm2} covers this observation.

%
\begin{proposition}
\label{proptoo-many-edges}
If a mixed graph $G=(V,D,B)$ with $V=[m]$ has $|D|+|B|>{{m}\choose{2}}$
edges, then $G$ is HTC-infinite-to-one.
\end{proposition}
\begin{pf}
Suppose $G$ is not HTC-infinite-to-one. Then there exists subsets
$(Y_v\dvtx v\in V)$, where each $Y_v$ satisfies the half-trek
criterion with respect to $v$ and for any pair of sets
$(Y_v,Y_w)$ it holds that $v\in Y_w$ implies $w\notin
Y_v$.

Fix a node $v\in V$. For every directed edge $u\to v\in D$, there
is a corresponding node $y\in Y_v$ for which it holds, by
Definition~\ref{defhalf-trek-crit}, that $y\bi v\notin B$.
Therefore, if there are $d_v$ directed edges pointing to $v$, then
there are $d_v$ nodes, namely, the ones in $Y_v$, that are not
adjacent to $v$ in the bidirected part $(V,B)$. If we consider
another node $w\in V$, with $d_w$ parents, then there are again
$d_w$ nonadjacencies $\{u,w\}$, $u\in Y_w$, in the bidirected
part. Moreover, $\{v,w\}$ cannot appear as a nonadjacency for both
node $v$ and node $w$ because of the requirement that $v\in Y_w$
imply $w\notin Y_v$. We conclude that there are at least $|D|$
nonedges in the bidirected part. In other words, $|D|+|B|\le
{{m}\choose{2}}$.~%
\end{pf}

We conclude the discussion of Theorems~\ref{mainthm} and
\ref{mainthm2} by pointing out that HTC-identifiability is equivalent
to a seemingly weaker criterion.\vadjust{\goodbreak}

%
\begin{definition}
\label{defweak-half-trek-crit}
A set of nodes $Y\subset V$ satisfies the weak half-trek criterion
with respect to node $v\in V$ if
\begin{longlist}[(iii)]
\item[(i)] $|Y|=|\pa(v)|$,
\item[(ii)] $Y\cap(\{v\}\cup\sib(v) ) = \varnothing$, and
\item[(iii)] there is a system of treks with no sided intersection
from $Y$ to $\pa(v)$ such that for any $w\in Y\cap\htr(v)$, the
trek originating at $w$ is a half-trek.
\end{longlist}
\end{definition}

%
\begin{lemma}\label{lemhalfornot}
Suppose the set $W\subset V$ satisfies the weak
half-trek criterion with respect to some node $v$. Then there
exists a set $Y$ satisfying the half-trek criterion with
respect to $v$, such that $Y\cap\htr(v)=W\cap
\htr(v)$.
\end{lemma}

Lemma~\ref{lemhalfornot} yields the following result; both the lemma
and the theorem are proved in Section~7 of the
supplement [\citet{Supplement}].

%
\begin{theorem}[(Weak HTC)]\label{mainthmweak}
Theorems~\ref{mainthm} and~\ref{mainthm2} hold when using the weak
half-trek criterion instead of the half-trek criterion. Moreover, a
graph $G$ can be proved to be rationally identifiable (or
generically infinite-to-one) using the weak half-trek criterion if
and only if $G$ is HTC-identifiable (or HTC-infinite-to-one).
\end{theorem}

\section{G-criterion}
\label{secG-crit}

The G-criterion, proposed in \citet{brito2006}, is a sufficient
criterion for rational identifiability in acyclic mixed graphs. The
criterion attempts to prove the fiber $\mathcal{F}(\Lambda,\Omega)$ to
be equal to $\{(\Lambda,\Omega)\}$ by solving the equation system
\[
\Sigma= (I-\Lambda)^{-T}\Omega(I-\Lambda)^{-1}
\]
in a stepwise manner. The steps yield the entries in $\Lambda$
column-by-column and, simultaneously, more and more rows and columns
for principal submatrices of $\Omega$. As explained in
Section~\ref{secproofs-H}, the half-trek method from
Section~\ref{secmain-results} starts from an equation system that has
$\Omega$ eliminated and then only proves~$\Lambda$ to
be uniquely identified. In this section, we show that, due to this
key simplification, the sufficient condition in the half-trek method
improves the G-criterion for acyclic mixed graphs.

To prepare for a comparison of the two criteria, we first restate
the identifiability theorem underlying the G-criterion in our own
notation. Enumerate the vertex set of an acyclic mixed graph $G$
according to any topological ordering as $V=[m]=\{1,\ldots,m\}$. (Then
$v\to w$ only if $v<w$.) Use the ordering to uniquely associate
bidirected edges to individual nodes by defining, for each $v \in V$,
the sets of siblings $S_<(v)=\{w\in\sib(v)\dvtx w<v\}$ and
$S_>(v)=\{w\in\sib(v)\dvtx w>v\}$. For a trek $\pi$, we write $t(\pi)$
to denote the target node, that is, $\pi$ is a trek from some node to
$t(\pi)$.

%
\begin{definition}[{[\citet{brito2006}]}]\label{defGcr}
A set of nodes $A\subset V$ satisfies the G-criterion with
respect to a node $v\in V$ if $A\subset V\setminus\{v\}$ and $A$ can
be partitioned into two (disjoint) sets $Y,Z$ with\vadjust{\goodbreak} $|Y|=|\pa(v)|$ and
$|Z|=|S_<(v)|$, with two systems of treks $\Pi\dvtx Y\rightrightarrows
\pa(v)$ and $\Psi\dvtx Z\rightrightarrows S_<(v)$, such that the following
condition holds:

If each trek $\pi\in\Pi$ is extended to a path
$\pi'$ by adding the edge $t(\pi)\to v$ to the right-hand side, and
each trek $\psi\in\Psi$ is similarly extended using $t(\psi)\bi v$,
then the set of paths $\{ \pi' \dvtx \pi\in\Pi\}\cup\{\psi' \dvtx \psi
\in\Psi\}$ is a set of
treks that has no sided intersection except at the common target
node $v$.
\end{definition}

Note that the paths $\pi'$ for $\pi\in\Pi$ are always treks. For
$\psi\in\Psi$, the requirement that $\psi'$ is a trek means that
$\psi$ cannot have an arrowhead at its target node.

For the statement of the main theorem about identifiability using the
G-criterion, define the depth of a node $v$ to be the length of the
longest directed path terminating at $v$. This number is denoted by
$\operatorname{Depth}(v)$.

%
\begin{theorem}[{[\citet{brito2006}]}]
\label{thmGcr}
Suppose $(A_v\dvtx v\in V)$ is a family of subsets of the vertex set
$V$ of an acyclic mixed graph $G$ and, for each $v$, the set~$A_v$ satisfies the G-criterion with respect to $v$. Then $G$
is rationally identifiable if at least one of the following two
conditions is satisfied:
\begin{longlist}[(C1)]
\item[(C1)] For all $v$ and all $w\in A_v$, it holds that
$\operatorname{Depth}(w)<\operatorname{Depth}(v)$.
\item[(C2)] For all $v$ and all $w\in A_v\cap( \htr(v) \cup
S_>(v) )$, the trek associated to node~$w$ in the definition
of the G-criterion is a half-trek. Furthermore, there is a total
ordering $\prec$ on $V$, such that if $w\in A_v\cap
( \htr(v)\cup S_>(v) )$, then $w\prec v$.
\end{longlist}
\end{theorem}

We remark that the ordering $\prec$ in condition (C2) need not agree
with any topological ordering of the graph. When using only condition
(C1) the theorem was given in \citet{brito2002g}, and the literature
is not always clear on which version of the G-criterion is concerned.
For instance, all examples in \citet{chan2010} can be proven to be
rationally identifiable by means of Theorem~\ref{thmGcr} as stated here.

We now compare the G-criterion to the half-trek criterion. We say
that a graph~$G$ is GC-identifiable if it satisfies the conditions of
Theorem~\ref{thmGcr}. The next theorem and proposition are proved in
Section~4 of the supplement [\citet{Supplement}].
They demonstrate that the half-trek method provides an improvement over
the G-criterion even for acylic mixed graphs.

%
\begin{theorem}\label{Gvshalf}
A GC-identifiable acyclic mixed graph is also HTC-identi\-fiable.
\end{theorem}

The graph in Figure~\ref{figtrekexamples} is HTC-identi\-fiable, as
was shown in Example~\ref{exGvshalfexample-HTCpart}.

%
\begin{proposition}
\label{propGvshalfexample}
The acyclic mixed graph in Figure~\ref{figtrekexamples} is not
GC-identi\-fiable.\vadjust{\goodbreak}
\end{proposition}


\begin{figure}

\includegraphics{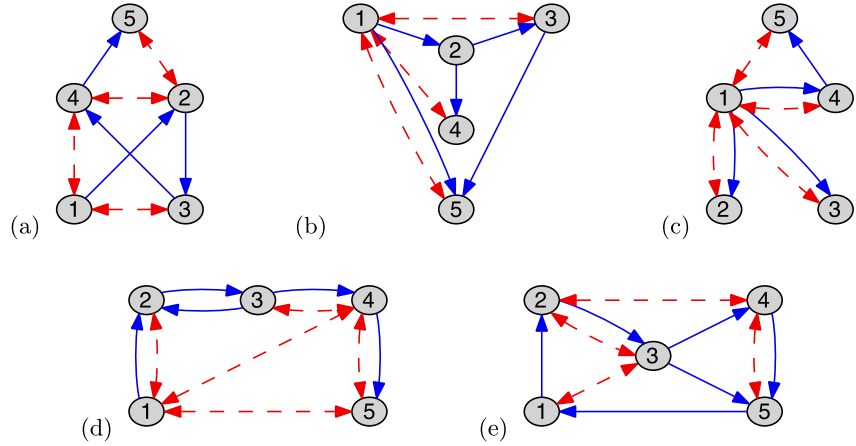}

\caption{Rationally identifiable mixed
graphs.}\label{figidexamples}
\end{figure}

\section{Examples}
\label{secexamples}

In the previous section the acyclic mixed graph from
Figure~\ref{figtrekexamples} was shown to be HTC-identifiable but
not GC-identifiable. In this section we give several other examples
that illustrate the conditions of our theorems and the ground that
lies beyond them. The examples are selected from the computational
experiments that we report on in Section~\ref{seccomputations}. We
begin with the identifiable class.

%
\begin{example}
Figure~\ref{figidexamples} shows 5 rationally identifiable mixed
graphs:
\begin{longlist}[(a)]
\item[(a)] This graph is simple and acyclic and, thus HTC- and
GC-identifiable; recall Proposition~\ref{propsimple-graphs}. There
are pairs $(\Lambda,\Omega)$ for which the fiber
$\mathcal{F}(\Lambda,\Omega)$ has positive dimension. By Theorem 2
in \citet{drton2011}, removing the edge $1\bi3$ would give a new
graph with all fibers of the form
$\mathcal{F}(\Lambda,\Omega)=\{(\Lambda,\Omega)\}$.
\item[(b)] The next graph is acyclic but not simple. It is HTC- and
GC-identifiable.
\item[(c)] This acyclic graph is HTC-inconclusive. The bidirected
part being connected, the example is not covered by the graph
decomposition technique discussed in
Section~\ref{secgraph-decompositions}.
\item[(d)] This is an example of a cyclic graph that is
HTC-identifiable.
\item[(e)] This cyclic graph is HTC-inconclusive.
\end{longlist}
\end{example}

On $m=5$ nodes, graphs with more than ${{5}\choose{2}}=10$ edges are
trivially generically infinite-to-one. The next example gives
nontrivial nonidentifiable graphs.

\begin{figure}

\includegraphics{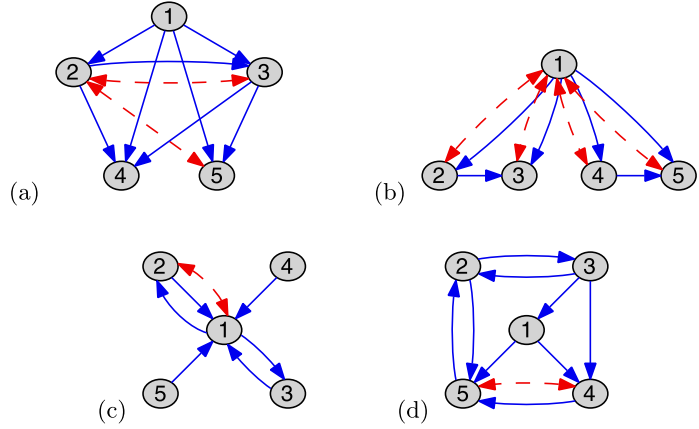}

\caption{Generically infinite-to-one
graphs.}\label{fignonidexamples}\vspace*{3pt}
\end{figure}

%
\begin{example}
All 4 graphs in Figure~\ref{fignonidexamples} are generically
infinite-to-one. The acyclic graph in (a) and the cyclic graph in
(c) are HTC-infinite-to-one. The acyclic graph in (b) and the
cyclic graph in (d) are HTC-inconclusive.\vadjust{\goodbreak}
\end{example}

\begin{figure}[b]
\vspace*{3pt}
\includegraphics{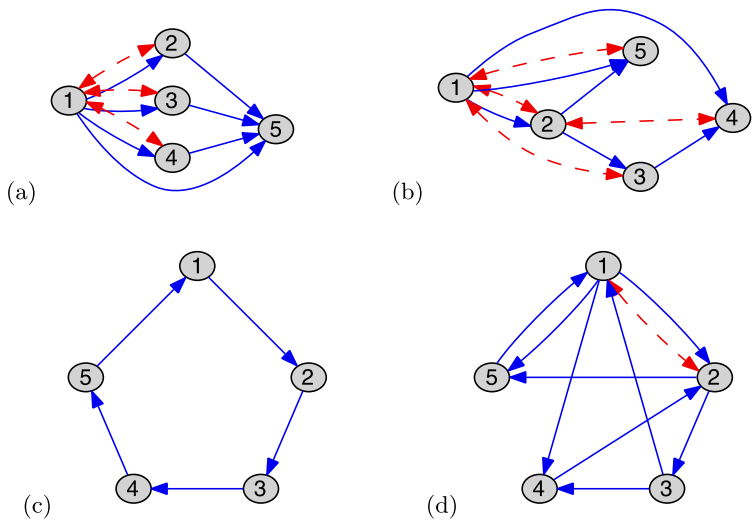}

\caption{Generically finite-to-one graphs.}\label{figfiniteidexamples}
\end{figure}

Many HTC-inconclusive graphs have fibers that are of
cardinality $2\le h<\infty$. An example of an acyclic 4-node graph
that is generically 2-to-one was given in \citet{brito2004}. Our next
example lists more graphs of this generically finite-to-one type.

%
\begin{example}
Figure~\ref{figfiniteidexamples} shows four mixed graphs that are
HTC-inconclu\-sive and not generically identifiable. All the graphs
have fibers that are generically finite:
\begin{longlist}[(a)]
\item[(a)] This graph is generically 2-to-1. We note that the
coefficients $\lambda_{v5}$, $v\in[4]$, can be identified; that is, any
two matrices $\Lambda,\Lambda'$ appearing in the same fiber have
an identical fifth column.
\item[(b)] Generically, the fibers of this graph have cardinality
of either one or three. For instance, let
\[
\omega_{11} = \cdots= \omega_{55} = 1,\qquad \omega_{12}
= \omega_{13} = \omega_{15} = \tfrac{1}{5},\qquad
\lambda_{23} =1.
\]
Define
\[
f(\lambda_{12}) = 529\lambda_{12}^4-460
\lambda_{12}^3- 3642\lambda_{12}^2-2380
\lambda_{12}-4271.
\]
Then, not considering the nongeneric situation with
$f(\lambda_{12})=0$, we have
\[
\bigl|\mathcal{F}(\Lambda,\Omega)\bigr| = \cases{ 3, &\quad$\mbox{if } f(\lambda_{12})
> 0,$\vspace*{2pt}
\cr
1, &\quad$\mbox{if } f(\lambda_{12}) < 0.$}
\]
The polynomial $f$ has two roots which are approximately $-2.16$ and
$3.44$.
\item[(c)] As shown in \citet{drton2011}, a cycle of length 3 or more
is generically 2-to-1.
\item[(d)] The next graph is not generically identifiable.
Generically, its fibers have at least two elements but not more than
10. Using the terminology from Definition~\ref{defid-degree}
below, the graph has degree of identifiability 10. We do not know of
an example of a fiber with more than two elements.
\end{longlist}
\end{example}


\section{Efficient algorithms for HTC-classification}
\label{secalgorithms}

While purely combinatorial, the identifiability conditions from
Theorems~\ref{mainthm} and~\ref{mainthm2} are not in a form that is
directly amenable to efficient computation. However, as we show in
this section, there exist polynomial-time algorithms for deciding
whether a~mixed graph $G$ is HTC-identifiable and whether $G$ is
HTC-infinite-to-one.
In the related context of the G-criterion, Chapter 4 in
\citet{brito2004} describes how the problem of determining the
existence of a set of nodes $Y$ satisfying the G-criterion with
respect to a given node $v$ can be solved by computation of maximum
flow in a derived directed graph. Our work for HTC-identifiability
extends this construction, which enables us to use maximum flow
computations to completely determine HTC-identifiability of
a mixed graph $G$. Furthermore, we show that whether $G$ is
HTC-infinite-to-one can be decided via a single max-flow
computation.

We first give some background on the max-flow problem; see
\citet{ford1962} and \citet{cormen2001}. Let $G=(V,D)$ be a directed
graph (or
``network'') with designated source and sink nodes $s,t\in V$. Let
$c_V\dvtx V\rightarrow\mathbb{R}_{\geq0}$ be a node-capacity function,
and let $c_D\dvtx D\rightarrow\mathbb{R}_{\geq0}$ be an edge-capacity
function. Then a \textit{flow} $f$ on $G$ is a function $f\dvtx
D\rightarrow
\mathbb{R}_{\geq0}$ that satisfies
\[
\sum_u f(u,v)=\sum_w
f(v,w)\leq c_V(v)
\]
for all nodes $v\neq s,t$, and
\[
f(u,v)\leq c_D(u,v)\vadjust{\goodbreak}
\]
for all edges $u\dir v\in D$.
The size $|f|$ of a flow $f$ on $G$ is the total amount of flow
passing from the source $s$ to the sink $t$, that is,
\[
|f|:=\sum_w f(s,w) = \sum
_u f(u,t) .
\]
The \textit{max-flow problem} on $(G,s,t,c_V,c_D)$ is the problem of
finding a flow $f$ whose size $|f|$ is maximum.

The computational complexity of the max-flow problem is known to be of
order $\mathcal{O}(|V|^3)$ if $G$ has no \textit{reciprocal} edge pairs.
A reciprocal edge pair consists of the two edges $v\dir w$ and $w\dir
v$ for distinct nodes $v\neq w$. (``Antiparallel'' is another term used
for such edge pairs.) In general, the complexity is
$\mathcal{O}((|V|+r)^3)$, where $r\leq|D|/2$ is the number of
reciprocal edge pairs. It is also known that if $c_V$ and $c_D$ are
both integer-valued, then there exists a maximal flow $f$ that is
integer-valued, and can be interpreted as a~sum of directed paths from
$s$ to $t$ with a flow of size~$1$ along each path
[\citet{ford1962}, \citet{cormen2001}]. (We note that the max-flow problem is
usually defined without bounded node capacities and on graphs with no
reciprocal edge pairs, but the more general problem stated here can be
converted to the standard form; see Section~6 of
the supplement [\citet{Supplement}] for details.)

\subsection{Deciding HTC-identifiability}

To determine whether a mixed graph $G=(V,D,B)$ is HTC-identifiable, we
first need to address the following subproblem. Given a node $v\in
V$, and a subset of ``allowed'' nodes $A\subseteq
V\setminus(\{v\}\cup\sib(v) )$, how can we efficiently
determine whether there exists a subset $Y\subseteq A$ satisfying the
half-trek criterion with respect to $v$? We now show that answering
this question is equivalent to solving a max-flow problem on a network
$G_{\mathrm{flow}}(v,A)$ with at most $2|V|+2$ nodes and at most
$3|V|+|D|+|B|$ edges.

\begin{figure}

\includegraphics{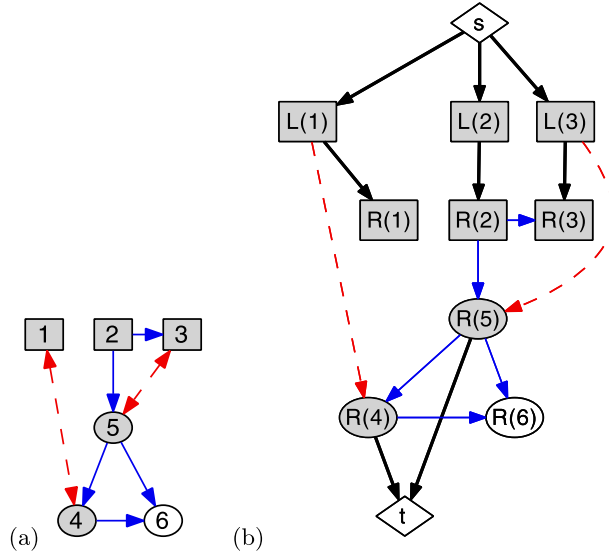}

\caption{Using max-flow to find a set satisfying the half-trek
criterion, for node $v=6$ and allowed nodes $A=\{1,2,3\}$. \textup{(a)}
The concerned mixed graph $G$. \textup{(b)} The network
$G_{\mathrm{flow}}(v,A)$.}\label{figMaxFlowHTC-id}
\end{figure}

We construct the network as follows; an example is shown in Figure
\ref{figMaxFlowHTC-id}. The vertex set of $G_{\mathrm{flow}}(v,A)$
comprises three types of nodes, namely,
\begin{longlist}[(a)]
\item[(a)] a source $s$ and a sink $t$,
\item[(b)] a ``left-hand copy'' $L(a)$ for each $a\in A$, and
\item[(c)] a ``right-hand copy'' $R(w)$ for each $w\in V$.
\end{longlist}
The edges of $G_{\mathrm{flow}}(v,A)$ are given by the following:
\begin{longlist}[(a)]
\item[(a)] $s\dir L(a)$ and $L(a)\dir R(a)$ for each $a\in A$ (thick
solid edges, in Figure~\ref{figMaxFlowHTC-id}),
\item[(b)] $L(a)\dir R(w)$ for each $a\bidir w\in B$ (dashed edges),
\item[(c)] $R(w)\dir R(u)$ for each $w\dir u\in D$ (solid edges), and
\item[(d)] $R(w)\dir t$ for each $w\in\pa(v)$ (thick solid edges).
\end{longlist}
Finally, we define the capacity functions. All edges have capacity
$\infty$. The source~$s$ and sink $t$ have capacity $\infty$, and all
other nodes have capacity $1$.\vadjust{\goodbreak}

The intuition for our construction is that a half-trek of the form
$y\dir x_1\dir\cdots\dir x_n=p$, with $y\in A$ and $p\in\pa(v)$,
will appear in the flow network as
\[
s\dir L(y)\dir R(y)\dir R(x_1)\dir\cdots\dir R(x_n)\dir
t ,
\]
and a half-trek of the form $y\bidir x_1\dir\cdots\dir x_n=p$ will
appear as
\[
s\dir L(y)\dir R(x_1)\dir\cdots\dir R(x_n)\dir t .
\]
By construction, no flow can exceed $|\pa(v)|$ in
size. Therefore, for practical purposes, all infinite capacities can
equivalently be replaced with capacity $|\pa(v)|$.

\begin{algorithm}
\caption{Testing HTC-identifiability of a mixed graph}
\label{algHTC-id}
\begin{algorithmic}
\STATE\textbf{Input:} $G=(V,D,B)$, a mixed graph on $m$ nodes\\
\STATE\textbf{Initialize:} $\texttt{SolvedNodes}\leftarrow\{v\dvtx\pa
(v)=\varnothing\}$.\\
\REPEAT
\FOR{$v=1,2,\ldots,m$}
\IF{$v\notin\texttt{SolvedNodes}$}
\STATE$A\leftarrow(\texttt{SolvedNodes}\cup(V\setminus\htr
(v)) )\setminus(\{v\}\cup\sib(v) )$.
\IF{$\texttt{MaxFlow}(G_{\mathrm{flow}}(v,A))=|\pa(v)|$}
\STATE$\texttt{SolvedNodes}\leftarrow\texttt{SolvedNodes}\cup\{v\}$.
\ENDIF
\ENDIF
\ENDFOR
\UNTIL{$\texttt{SolvedNodes}=V$ or no change has occurred in the last
iteration.}
\STATE\textbf{Output:} ``yes'' if $\texttt{SolvedNodes}=V$, ``no''
otherwise.
\end{algorithmic}
\end{algorithm}

The following theorem is proved in Section~6
of the supplement.

%
\begin{theorem}\label{thmMaxFlowHTC-id}
Given a mixed graph $G=(V,D,B)$, a~node $v\in V$ and a~subset of
``allowed'' nodes $A\subseteq V\setminus(\{v\}\cup
\sib(v) )$, there exists a set $Y\subseteq A$ satisfying the
half-trek criterion with respect to $v$ if and only if the flow
network $G_{\mathrm{flow}}(v,A)$ has maximum flow equal to
$|\pa(v)|$.
\end{theorem}

Using Theorem~\ref{thmMaxFlowHTC-id}, we are able to
give an algorithm to determine whether $G$ is HTC-identifiable. If $G$
is HTC-identifiable, then, by Definition~\ref{defhalf-trek-crit}, we
have an ordering $\prec$ on $V$, and for each $v$, a set $Y_v$
satisfying the half-trek criterion with respect to $v$, such that any
$w\in Y_v\cap\htr(v)$ must be $\prec v$. Therefore, by Theorem
\ref{thmMaxFlowHTC-id}, the network $G_{\mathrm{flow}}(v,A)$ must
have maximum flow size $|\pa(v)|$, where $A$\vadjust{\goodbreak} is the set of nodes that
are ``allowed'' to be in $Y_v$ according to the ordering $\prec$, that is,
\[
A= \bigl[\{w\dvtx w\prec v\} \cup\bigl(V\setminus\htr(v)\bigr) \bigr]
\setminus
\bigl[\{v\}\cup\sib(v) \bigr] .
\]
This intuition is formalized in Algorithm~\ref{algHTC-id}. In
Section~6 of the supplement, we prove the
following theorem, which
states that Algorithm~\ref{algHTC-id} correctly determines
HTC-identifiability.

%
\begin{theorem}
\label{thmMaxFlowAlgHTC-id}
A mixed graph $G=(V,D,B)$ is HTC-identifiable if and only if
Algorithm~\ref{algHTC-id} returns ``yes.'' Furthermore, the
algorithm has complexity at most $\mathcal{O}(|V|^2(|V|+r)^3)$,
where $r\leq|D|/2$ is the number of reciprocal edge pairs in $D$.
\end{theorem}

\subsection{Deciding if a graph is HTC-infinite-to-one}
\label{subsecGflow}

To determine whether a mixed graph $G=(V,D,B)$ is HTC-infinite-to-one,
we may again appeal to max-flow computation. It now suffices to solve
a single larger max-flow problem, with at most $\frac{3}{2} |V|^2+2$
nodes and at most $|V|\cdot(\frac{3}{2} |V|+2|D|+|B|)$ edges, and
$|V|\cdot r$ reciprocal edge pairs, where $r$ is the number of
reciprocal edge pairs in $G$.

\begin{figure}

\includegraphics{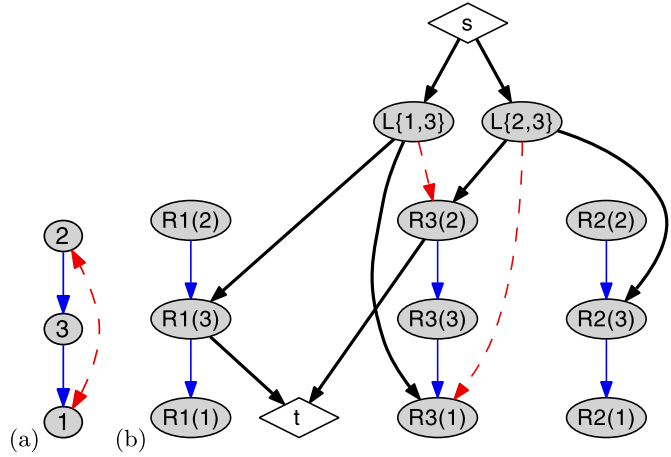}

\caption{Using max-flow to test whether a mixed graph is
HTC-infinite-to-one. \textup{(a)}~A~mixed graph $G$ on 3 nodes. \textup
{(b)} The
associated flow network $G_{\mathrm{flow}}$.}\label{figMaxFlowHTC-nonid}
\end{figure}

The relevant flow network $G_{\mathrm{flow}}$ is constructed as
follows; an example is shown in Figure~\ref{figMaxFlowHTC-nonid}.
The nodes of $G_{\mathrm{flow}}$ are as follows:
\begin{longlist}[(a)]
\item[(a)] a source $s$ and a sink $t$,
\item[(b)] a ``left-hand copy'' $L\{v,w\}$ for each unordered pair
$\{v,w\}\subset V$ with $v\bidir w\notin B$, and
\item[(c)] a ``right-hand copy'' $R_v(w)$ for each $v,w\in V$.
\end{longlist}
The edges of $G_{\mathrm{flow}}$ are as follows:
\begin{longlist}[(a)]
\item[(a)] $s\dir L\{v,w\}$ and $L\{v,w\}\dir R_v(w)$ for each unordered
pair $\{v,w\}\subset V$ with $v\bidir w\notin B$ (thick solid edges, in
Figure~\ref{figMaxFlowHTC-nonid}),\vadjust{\goodbreak}
\item[(b)] $L\{v,w\}\dir R_v(u)$ for each $v,w,u$ with $v\neq w$ such
that $v\bidir w\notin B$ but $w\bidir u\in B$ (dashed edges),
\item[(c)] $R_v(w)\dir R_v(u)$ for each $v,w,u\in V$ with $w\dir u \in
D$ (solid
edges), and
\item[(d)] $R_v(w)\dir t$ for each $v,w\in V$ with $w\in\pa(v)$
(thick
solid edges).
\end{longlist}
Finally, the edge capacity function assigns capacity $\infty$ to all
edges, and the node capacity function gives capacity $\infty$ to the
source $s$ and sink $t$ and capacity $1$ to all other nodes. If
useful in practice, the infinite capacities can be set to $|V|^2$, as
no flow can have size larger than $|V|^2$.

The intuition for the construction just given is as follows. If the
mixed graph $G$ is not HTC-infinite-to-one, then simultaneously for
all nodes $v\in V$, we can find systems of half-treks with no sided
intersection $Y_v \rightrightarrows\pa(v)$, such that $Y_v$ does not
contain $v$ or any siblings of $v$, and $w\in Y_v$ implies \mbox{$v\notin
Y_w$}. Writing $y_{(v,k)}\eqbi z_{(v,k),1}$ to represent
either $y_{(v,k)}= z_{(v,k),1}$ or $y_{(v,k)}\bidir z_{(v,k),1}$, a half-trek
\[
\pi_{(v,k)}\dvtx y_{(v,k)}\eqbi z_{(v,k),1}\dir
z_{(v,k),2}\dir\cdots\dir k
\]
with $k\in\pa(v)$ and $y_{(v,k)}\in Y_v$ corresponds to a path in the
network $G_{\mathrm{flow}}$ given by
\[
\tilde{\pi}_{(v,k)}\dvtx s\dir L\{v,y_{(v,k)}\}\dir
R_v(z_{(v,k),1})\dir R_v(z_{(v,k),2})\dir
\cdots\dir R_v(k)\dir t .
\]
Therefore, in the maximum flow on $G_{\mathrm{flow}}$, if $\{v,w\}$ is
used by one of the paths passing through the $R_v(\cdot)$ copy of the
graph, then it will not get used by any of the flows passing through
the $R_w(\cdot)$ copy of the graph.

The following theorem is proved in Section~6
of the supplement.

%
\begin{theorem}\label{thmMaxFlowHTC-nonid}
A mixed graph $G=(V,D,B)$ is HTC-infinite-to-one if and only if
$G_{\mathrm{flow}}$ has maximum flow size strictly less than $|D|=\sum
_{v\in
V}|\pa(v)|$. The computational complexity of solving this
max-flow problem is\break $\mathcal{O}(|V|^3(|V|+r)^3)$, where $r\leq
|D|/2$ is the number of reciprocal edge pairs in $G$.
\end{theorem}


\section{Computational experiments}
\label{seccomputations}

This section reports on the results of an exhaustive study of all
mixed graphs with $m\le5$ nodes, for which the identification problem
can be fully solved by means of algebraic techniques. Moreover, we
show simulations in which we apply our new combinatorial criteria to
graphs with $m=25$ and $50$ nodes.

\subsection{Exhaustive computations on small graphs}

We applied the half-trek and the G-criterion as well as algebraic
techniques to all mixed graphs on $m\leq5$ nodes. All algebraic
computations were done with the software
\textsc{Singular} [\citet{Singular}]; see Section~1
of the supplement [\citet{Supplement}] for
details. The G-criterion and the max-flow
algorithms from Section~\ref{secalgorithms} were implemented in
\textsc{R}~[\citet{R}] and \textsc{MATLAB} [\citet{MATLAB2010}],
respectively.

%
\begin{table}[b]
\caption{Classification of unlabeled mixed graphs with $3\le m\le5$
nodes; column~``HTC''~gives~counts of
HTC-classifiable graphs}\label{Table345nodes}
\begin{tabular*}{\textwidth}{@{\extracolsep{\fill}}lcccccc@{}}
\hline
&\multicolumn{2}{c}{$\bolds{m=3}$}&\multicolumn{2}{c}{$\bolds{m=4}$}&
\multicolumn{2}{c@{}}{$\bolds{m=5}$}\\[-6pt]
&\multicolumn{2}{c}{\hrulefill}&\multicolumn{2}{c}{\hrulefill}&
\multicolumn{2}{c@{}}{\hrulefill}\\
\textbf{Unlabeled mixed graphs}&\textbf{Total}&\textbf{HTC}&\textbf
{Total}&\textbf{HTC}&\textbf{Total}&\textbf{HTC}\\
\hline
Acyclic, $\leq{m\choose2}$ edges&\textit{22}&& \phantom{0}\textit
{715}&& \phantom{0,}\textit{103,670}&\\
\quad rationally identifiable& 17&17& \phantom{0}343&343& \phantom
{00,}32,378&\phantom{0}32,257\\
\quad generically finite-to-one& \phantom{0}0&--& \phantom{000}4&--& \phantom
{000.,}1166&--\\
\quad generically $\infty$-to-one& \phantom{0}5& \phantom{0}5& \phantom{0}368
&368& \phantom{00,}70,126& \phantom{0}70,099\\[3pt]
Acyclic, $>{m\choose2}$
edges&\textit{18}&&\phantom{0}\textit{852}&&\phantom{0,}\textit
{152,520}&\\[6pt]
Cyclic, $\leq{m\choose2}$ edges&\phantom{0}\textit{6}&&\phantom
{0}\textit{718}&&\phantom{0,}\textit{348,175}&\\
\quad rationally identifiable& \phantom{0}2&\phantom{0}2& \phantom{0}239&230&
\phantom{00,}91,040&\phantom{0}78,586\\
\quad generically finite-to-one& \phantom{0}1&--& \phantom{00}75&---& \phantom
{00,}44,703&--\\
\quad generically $\infty$-to-one& \phantom{0}3& \phantom{0}3& \phantom
{0}404&383& \phantom{0,}212,432&
202,697\\[3pt]
Cyclic, $>{m\choose2}$ edges&\textit{58}&&\textit{9307}&&\textit
{8,439,859}&\\
\hline
\end{tabular*}
\end{table}

The results are given in Table~\ref{Table345nodes}, where we treat
graphs as unlabeled, that is, we count isomorphism classes of graphs
with respect to permutation of the vertex set $V=[m]$. The table
distinguishes between acyclic and cyclic (i.e., nonacyclic)
graphs. In each case, we single out the graphs with more than
${{m}\choose{2}}$ edges. These are trivially generically infinite-to-one
and also HTC-infinite-to-one according to
Proposition~\ref{proptoo-many-edges}. The remaining graphs are
classified into three disjoint groups, namely, rationally identifiable
graphs, generically infinite-to-one graphs and generically
finite-to-one graphs. The following notion makes the distinctions and
terminology precise. Here, $\mathbb{C}^D_\mathrm{reg}$ is
defined as $\mathbb{R}^D_\mathrm{reg}$ but allowing for complex matrix
entries. We write $\mathbb{C}^{m\times m}_\mathrm{sym}$ for the space
of symmetric $m\times m$ complex matrices.

%
\begin{definition}
\label{defid-degree}
Let $G=(V,D,B)$ be a mixed graph. Then the complex rational map
$\phi_{G,\mathbb{C}}$, obtained by extending the map $\phi_G$ to
$\mathbb{C}^D_\mathrm{reg}\times\mathbb{C}^{m\times
m}_\mathrm{sym}$, is generically $h$-to-one with
$h\in\mathbb{N}\cup\{\infty\}$, and we call $h=\operatorname{ID}(G)$ the
degree of identifiability of~$G$.
\end{definition}

A mixed graph $G$ is rationally identifiable if and only if its degree
of identifiability $\operatorname{ID}(G)=1$. Similarly, $G$ is generically
infinite-to-one if and only if $\operatorname{ID}(G)=\infty$; in that case
the fiber $\mathcal{F}(\Lambda,\Omega)\subset
\mathbb{R}^D_\mathrm{reg}\times\operatorname{PD}(B)$ defined in
(\ref{eqfiber}) is generically of positive dimension. In
Table~\ref{Table345nodes}, a graph $G$ is generically finite-to-one
if $2\le\operatorname{ID}(G) <\infty$ and, thus,
$\mathcal{F}(\Lambda,\Omega)$ is generically finite with
$|\mathcal{F}(\Lambda,\Omega)|\le\operatorname{ID}(G)$. If $\operatorname{ID}(G)$
is finite and even, $G$ cannot be generically identifiable
because polynomial equations have complex solutions appearing in
conjugate pairs and $\mathcal{F}(\Lambda,\Omega)$ always contains at
least one (real) point, that is,\  $(\Lambda,\Omega)$.
If $\operatorname{ID}(G)$ is odd, we cannot exclude the possibility
that the equation defining $\mathcal{F}(\Lambda,\Omega)$
generically only has one real point, leading to generic
identifiability. However, we did not observe this in any examples we
checked.\looseness=1

Table~\ref{Table345nodes} shows that our half-trek method yields a
perfect classification of acyclic graphs with $m\le4$ nodes and
cyclic graphs with $m\le3$ nodes. Among the acyclic graphs with
$m=5$ nodes, our method misses 121 rationally identifiable graphs and
27 generically infinite-to-one graphs. The gaps are larger for cyclic
graphs, but the method still classifies 86\% of
the rationally identifiable graphs correctly and misses less than 5\%
of the generically infinite-to-one graphs. In the supplementary
article
[\citet{Supplement}],
we list some rationally identifiable graphs and some generically infinite-to-one
graphs that are not classifiable using our method (i.e.,\ that are
HTC-inconclusive). The degree of
identifiability~$\operatorname{ID}(G)$ of a graph~$G$ with $5$ nodes can be
any number in $[8]\cup\{10\}$, and any number in~$[4]$ when $G$ is
acyclic. For example, the graphs in
Figure~\ref{figfiniteidexamples}(a), (b) and (d) have
$\operatorname{ID}(G)=2$, $3$ and $10$, respectively.

We also tracked which acyclic graphs are rationally
identifiable according to the G-criterion from Theorem~\ref{thmGcr}.
Since this method depends on the choice of a topological ordering of
the nodes, we tested each possible topological ordering.
Our computation shows that the G-criterion finds all rationally
identifiable acyclic graphs with $m\le4$ nodes. For $m=5$, the
G-criterion proves 31,830 acyclic graphs to be rationally
identifiable but misses 427 of the HTC-identifiable acyclic
graphs.

\subsection{Simulations for large graphs}

Exhaustive computations become prohibitive for more than 5 nodes.
Furthermore, algebraic computations are not feasible for larger
graphs. Instead, we test the HTC-status of randomly generated mixed
graphs with $m=25$ or $m=50$ nodes.

\begin{figure}

\includegraphics{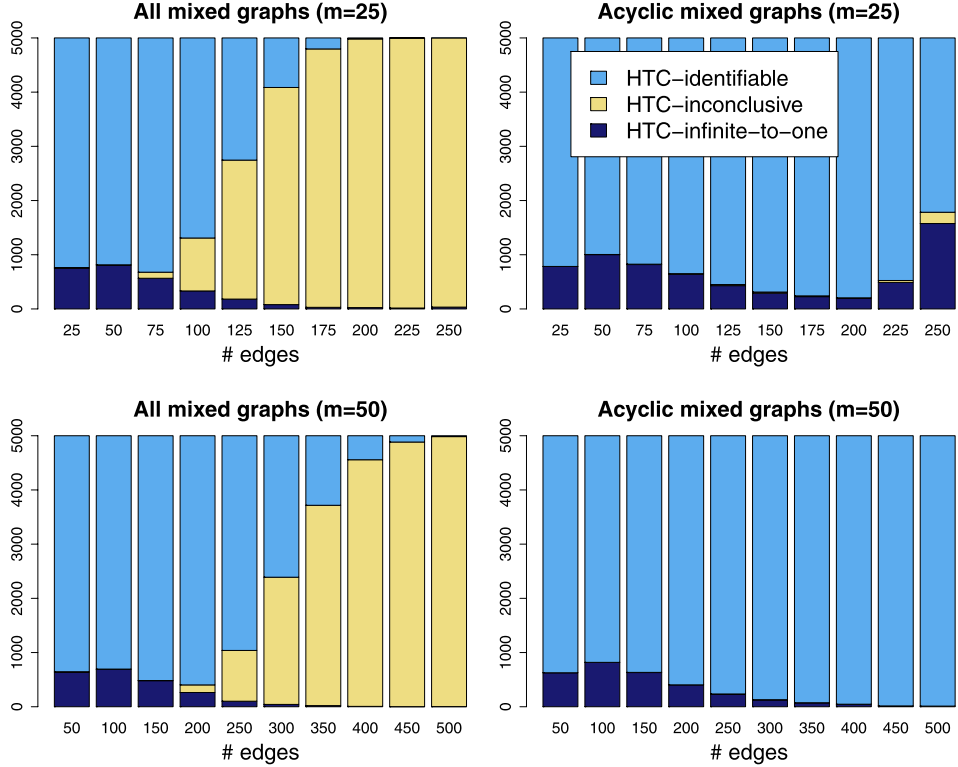}

\caption{Classification of labeled mixed graphs
with $m=25$ and $m=50$ nodes. Each bar represents 5000 randomly
drawn graphs
with fixed number of edges, ranging from $m$ to $m\cdot10$.}
\label{Figurem25m50}
\end{figure}

For each value $n=k\cdot m$ for $k\in[10]$, we randomly sampled $5000$
labeled mixed graphs on $m$ nodes with $n$ edges, by selecting a
subset of size $n$ from the set of all possible edges, which consists
of $2\cdot{{m}\choose{2}}$ directed edges and ${m\choose2}$ bidirected
edges. We repeated this process with acyclic graphs only; the choice
is then from ${m\choose2}$ directed edges and ${m\choose2}$
bidirected edges. The results of these simulations are shown in
Figure~\ref{Figurem25m50}. When the graphs are restricted to be
acyclic, most are HTC-identifiable and only extremely few are
HTC-inconclusive. When we do not restrict to acyclic graphs, on the
other hand, we see that as the number of edges increases, the
proportion of HTC-inconclusive graphs grows rapidly.

\begin{figure}

\includegraphics{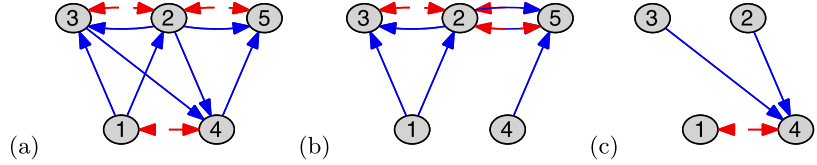}

\caption{An acyclic mixed graph shown in \textup{(a)} and its two mixed
components shown in \textup{(b)} and \textup{(c)}.}
\label{fighalf-trek-decomposition}
\end{figure}

\section{Decomposition of acyclic graphs}
\label{secgraph-decompositions}

In this section we discuss how, for acyclic graphs, the scope of
applicability of our half-trek method can be extended via a graph
decomposition due to \citet{tian2005}. Let $G=(V,D,B)$ be an acyclic
mixed graph, and let $C_1,\ldots,C_k\subset V$ be the (pairwise
disjoint) vertex sets of the connected components of the bidirected
part $(V,B)$. For $j\in[k]$, let $B_j = B\cap(C_j\times C_j)$ be the
bidirected edges in the $j$th connected component. Let $V_j$ be
the union of $C_j$ and any parents of nodes in $C_j$, that is,
\[
V_j = C_j \cup\bigl\{ \pa(v) \dvtx v\in
C_j \bigr\},\qquad j=1,\ldots,k.
\]
Clearly, the sets $V_1,\ldots,V_k$ need not be pairwise disjoint. Let
$D_j$ be the set of edges $v\to w$ in the directed part $(V,D)$ that
have $v\in V_j$ and $w\in C_j$. The decomposition of \citet{tian2005}
involves the graphs $G_j=(V_j,D_j,B_j)$, for $j\in[k]$. We refer to
these as the \textit{mixed components} $G_1,\ldots,G_k$ of $G$.
Figure~\ref{fighalf-trek-decomposition} gives an example.

The mixed components $G_1,\ldots,G_k$ create a partition of the edges
of $G$. There is an associated partition of the entries of
$\Lambda\in\mathbb{R}^D$ that yields submatrices
$\Lambda_1,\ldots,\Lambda_k$ with each $\Lambda_j\in\mathbb{R}^{D_j}$;
recall that for an acyclic graph
$\mathbb{R}^D_\mathrm{reg}=\mathbb{R}^D$. Similarly, from
$\Omega\in\operatorname{PD}(B)$, we create matrices
$\Omega_1,\ldots,\Omega_k$ with each $\Omega_j\in\operatorname{PD}(B_j)$,
where $\operatorname{PD}(B_j)$ is defined with respect to the graph $G_j$,
that is, the set contains matrices indexed by $V_j\times V_j$. We
define $\Omega_j$ by taking the submatrix $\Omega_{C_j,C_j}$ from
$\Omega$ and extending it by setting $(\Omega_j)_{vv}=1$ for all
$v\in
V_j\setminus C_j$. The work leading up to Theorems 1 and 2 in
\citet{tian2005} shows that, for all $j\in[k]$, there is a rational
map $f_j$ defined on the entire cone of $m\times m$ positive definite
matrices such that
\[
f_j\circ\phi_G(\Lambda,\Omega) = \phi_{G_j}(
\Lambda_j,\Omega_j)
\]
for all $\Lambda\in\mathbb{R}^D$ and $\Omega\in\operatorname{PD}(B)$. In
turn, there is a rational map $g$ defined everywhere on the product of
the relevant cones of positive definite matrices such that
\[
g \bigl(\phi_{G_1}(\Lambda_1,\Omega_1),\ldots,
\phi_{G_k}(\Lambda_k,\Omega_k)\bigr) =
\phi_G(\Lambda,\Omega)
\]
for all $\Lambda\in\mathbb{R}^D$ and $\Omega\in\operatorname{PD}(B)$. We
thus obtain the following theorem.

%
\begin{theorem}
For an acyclic mixed graph $G$ with mixed components
$G_1,\ldots, G_k$, the following holds:
\begin{longlist}[(iii)]
\item[(i)] $G$ is rationally (or generically) identifiable if and
only if all components $G_1,\ldots,G_k$ are rationally (or generically)
identifiable;
\item[(ii)] $G$ is generically infinite-to-one if and only if there
exists a component~$G_j$ that is generically infinite-to-one;
\item[(iii)] if each $G_j$ is generically $h_j$-to-one with
$h_j<\infty$,
then $G$ is generically $h$-to-one with $h=\prod_{j=1}^k h_j$.
\end{longlist}
\end{theorem}

We remark that this theorem could also be stated as
$\operatorname{ID}(G)=\prod_{j=1}^k \operatorname{ID}(G_j)$, in terms of the
degree of identifiability from Definition~\ref{defid-degree}.

The next theorem makes the observation that when applying our
half-trek method to an acyclic graph, we may always first decompose
the graph into its mixed components, which may result into
computational savings.

%
\begin{theorem}
\label{thmHTC-decomposition}
If an acyclic mixed graph $G$ is HTC-identifiable, then all its mixed
components $G_1,\ldots,G_k$ are HTC-identifiable. Furthermore, $G$
is HTC-infinite-to-one if and only if there exists a mixed component $G_j$
that is HTC-infinite-to-one.
\end{theorem}
\begin{pf}
The claim about HTC-identifiability follows from
Lemma~4 in Section~5 of the supplement
[\citet{Supplement}].
The second statement is a consequence of Lemmas~5 and~6 from the same section.
\end{pf}

The benefit of graph decomposition goes beyond computation in that
some identification methods apply to all mixed components but not to
the original graph. In \citet{tian2005}, this is exemplified for the
G-criterion. More precisely, the 4-node example given there concerns
the early version of the G-criterion from \citet{brito2002g} that
includes only condition (C1) from Theorem~\ref{thmGcr} but not
condition (C2), which is due to \citet{brito2006}. However, graph
decomposition allows one to also extend the scope of our more general
half-trek method, where passing to mixed components can avoid problems
with finding a suitable total ordering of the vertex set.
Surprisingly, however, the extension is possible only for the
sufficient condition, that is, HTC-identifiability;
Theorem~\ref{thmHTC-decomposition} gives an equivalence result for
HTC-infinite-to-one graphs.

%
\begin{proposition}
\label{prophalf-trek-extend-decomposition}
The acyclic mixed graph in
Figure~\ref{fighalf-trek-decomposition}\textup{(a)} is not HTC-identifiable
but both its mixed components are HTC-identifiable.
\end{proposition}
\begin{pf}
Suppose for a contradiction that the original graph $G$ is
HTC-identifiable and that the sets $Y_3$, $Y_4$ and
$Y_5$ are part of the family of sets appearing in
Theorem~\ref{mainthm}. In particular, each set has two elements and
satisfies\vadjust{\goodbreak} the half-trek criterion with respect to its subscript.
Now, the presence of the edge $2\bi3$ implies that $Y_3\subset
\{1,4,5\}$. Moreover, $Y_3\not=\{1,4\}$ because the sole
half-trek from $4$ to $3$ has $1$ in its right-hand side and all
half-treks from $1$ to $3$ are directed paths and thus have the
source $1$ on their right-hand side as well. It follows that $5\in
Y_3$ and, thus, $3\notin Y_5$. Since $2\bi5$ is in $G$,
it must hold that $Y_5=\{1,4\}$. Examining the descendant sets
$ \htr(v)$, we see that the total ordering $\prec$ in
Theorem~\ref{mainthm} ought to satisfy $4\prec5\prec3$. Since
$1\in\sib(4)$ and $3,5\in\htr(4)$, we
conclude that $Y_4\subset\{2\}$, which is a contradiction
because $Y_4$ must have two elements.

Turning to the mixed components of $G$, it is clear that the
component shown in Figure~\ref{fighalf-trek-decomposition}(c) is
HTC-identifiable because it is a simple graph; recall
Proposition~\ref{propsimple-graphs}. The component in
Figure~\ref{fighalf-trek-decomposition}(b) is HTC-identifiable
because Theorem~\ref{mainthm} applies with the choice of
\[
Y_1 = Y_4=\varnothing,\qquad Y_2= \{1\},\qquad
Y_5 = \{1,4\},\qquad Y_3 = \{1,5\},
\]
and any ordering that respects $5\prec3$.
\end{pf}

As seen in Table~\ref{Table345nodes}, the half-trek method misses 121
rationally identifiable acyclic graphs with 5 nodes, among them is the
example from Proposition~\ref{prophalf-trek-extend-decomposition}.
After graph decomposition, the half-trek method proves 9 of the 121
examples to be rationally identifiable. The remaining 112 graphs all
have a~connected bidirected part; see Figure~\ref{figidexamples}(c)
for an example. On 5 nodes, there are 27 generically infinite-to-one
graphs that are HTC-inconclusive. All of these have a connected
bidirected part.

\section{Proofs for the half-trek criterion}
\label{secproofs-H}

In this section we prove the two main theorems stated in
Section~\ref{secmain-results}. We begin with the identifiability theorem.

%
\begin{reptheorem}[(HTC-identifiability)]
Let $(Y_v\dvtx v\in V)$ be a family of subsets of the vertex
set $V$ of a mixed graph $G$. If, for each node $v$, the set
$Y_v$ satisfies the half-trek criterion with respect to $v$, and
there is a total ordering $\prec$ on the vertex set $V$ such that
$w\prec v$ whenever $w\in Y_v\cap\htr(v)$, then $G$ is
rationally identifiable.
\end{reptheorem}
\begin{pf}
Let $\Sigma=\phi_G(\Lambda_0,\Omega_0)$ be a~matrix in the image of~$\phi_G$, given by a~generically chosen pair
$(\Lambda_0,\Omega_0)\in\Theta=\mathbb{R}^D_\mathrm{reg}\times
\operatorname{PD}(B)$. For generic identifiability, we need to show that
the equation
%
\begin{equation}
\label{eqfiber-eqn} \Sigma=(I-\Lambda)^{-T}\Omega(I-
\Lambda)^{-1}
\end{equation}
has a unique solution in $\Theta$, namely,
$(\Lambda,\Omega)=(\Lambda_0,\Omega_0)$. However, a pair
$(\Lambda,\Omega)$ solves~(\ref{eqfiber-eqn}) if and only if
%
\begin{equation}
\label{eqnzeros} \bigl[(I-\L)^T \Sigma(I-\L) \bigr]_{vw}=0\qquad
\forall(v,w)\notin B \mbox{ and } v\not= w
\end{equation}
and
%
\begin{equation}
\label{eqnnonzeros} \bigl[(I-\L)^T \Sigma(I-\L)
\bigr]_{vw}=\Omega_{vw}\qquad \forall(v,w)\in B \mbox{ or } v= w.
\end{equation}
The nonzero entries of $\Omega$ appearing in~(\ref{eqnnonzeros})
are freely varying real numbers that are subject only to the
requirement that $\Omega$ be positive definite. For cyclic graphs,
(\ref{eqfiber-eqn}) contains rational equations. Hence, the focus is
on~(\ref{eqnzeros}), which defines a polynomial equation system
even when the graph is cyclic.

We prove the theorem by solving the equations~(\ref{eqnzeros}) in
a stepwise manner according to the ordering $\prec$. When visiting
node $v$, the goal is to recover the $v$th column of $\Lambda$ as a
function of $\Sigma$. Based on solving linear equation systems, the
functions of $\Sigma$ that give the entries of $\Lambda$ will always
be rational functions, proving our stronger claim of rational (as
opposed to mere generic) identifiability.

For our proof we proceed by induction and assume that, for all
$w\prec v$, we have recovered the entries of the vector
$\L_{\pa(w),w}$ as (rational) expressions in~$\Sigma$. To solve
for $\L_{\pa(v),v}$, let $Y_v=\{y_1,\ldots,y_n\}$ and
$\pa(v)=\{p_1,\ldots,p_n\}$. Define $\AA\in\R^{n\times n}$ as
\[
\AA_{ij}=\cases{ \bigl[(I-\L)^T\Sigma
\bigr]_{y_ip_j},&\quad$\mbox{if } y_i\in\htr(v),$\vspace*{2pt}
\cr
\Sigma_{y_ip_j},&\quad$\mbox{if } y_i\notin\htr(v).$}
\]
Define $\bb\in\R^n$ as
\[
\bb_i=\cases{ \bigl[(I-\L)^T\Sigma
\bigr]_{y_iv},&\quad${\mbox{if } y_i\in\htr(v),}$\vspace*{2pt}
\cr
\Sigma_{y_iv},&\quad ${\mbox{if } y_i\notin\htr(v).}$}
\]
Note that both $\AA$ and $\bb$ depend only on $\Sigma$ and the
columns $\Lambda_{\pa(w),w}$ with $w\in Y_v\cap\htr(v)$, which
are assumed already to be known as a function of $\Sigma$ because
$w\in Y_v\cap\htr(v)$ implies $w\prec v$. We now claim that
the vector $\Lambda_{\pa(v),v}$ solves the equation system
$\AA\cdot\Lambda_{\pa(v),v}=\bb$.

First, consider an index $i$ with $y_i\in Y_v\cap\htr(v)$.
Since $Y_v$ satisfies the half-trek criterion with respect to
$v$, the node $y_i\not=v$ is not a sibling of $v$. Therefore, by
(\ref{eqnzeros}),
\[
\bigl[(I-\L)^T\ \Sigma(I-\L) \bigr]_{y_iv} =0
\quad\Longrightarrow\qquad\bigl[(I-\L)^T\Sigma\Lambda\bigr]_{y_iv}=
\bigl[(I-\L)^T\Sigma\bigr]_{y_iv}.
\]
It follows that
\begin{eqnarray*}
(\AA\cdot\Lambda_{\pa(v),v} )_i&= &\sum
_{j=1}^n \bigl[(I-\L)^T\Sigma)
\bigr]_{y_ip_j}\Lambda_{p_jv}\\
& =&
\bigl[(I-\L)^T\Sigma\Lambda\bigr]_{y_iv}= \bigl[(I-
\L)^T\Sigma\bigr]_{y_iv}=\bb_i.
\end{eqnarray*}

Second, let $i$ be an index with $y_i\in Y_v\setminus\htr(v)$. Then
\[
(\AA\cdot\Lambda_{\pa(v),v} )_i= \sum
_{j=1}^n \S_{y_ip_j}\Lambda_{p_jv} = [
\S\L]_{y_iv} = \bigl[(I-\L)^{-T}\O(I-\L)^{-1}\L
\bigr]_{y_iv}.
\]
By definition of $ \htr(v)$, we know that $[(I-\L)^{-T}\O]_{y_iv}=0$.
Adding this zero and using that $(I-\L)^{-1}=I+(I-\L)^{-1}\L$, we
obtain that
\begin{eqnarray*}
(\AA\cdot\Lambda_{\pa(v),v} )_i&=& \bigl[(I-\L)^{-T}
\O(I-\L)^{-1}\L\bigr]_{y_iv}+ \bigl[(I-\L)^{-T}\O
\bigr]_{y_iv}\\
&=&
\bigl[(I-\L)^{-T}\O(I-\L)^{-1} \bigr]_{y_iv} =
\S_{y_iv}=\bb_i.
\end{eqnarray*}
Therefore, $\AA\cdot\Lambda_{\pa(v),v}=\bb$, as claimed.

By Lemma~\ref{leminvert} below, the matrix $\AA$ is invertible in
the generic situation. Therefore, we have shown that
$\Lambda_{\pa(v),v}=\AA^{-1}\bb$ is a rational function of~$\Sigma$.
Proceeding inductively according to the vertex ordering
$\prec$, we recover~$\Lambda_{\pa(v),v}$ for all $v$ and, thus, the
entire matrix $\Lambda$, as desired.
\end{pf}

%
\begin{lemma}
\label{leminvert}
Let $v\in V$ be any node. Let $Y\subset V\setminus(\{v\}\cup
\sib(v) )$, with $|Y|=|\pa(v)|=n$. Write
$Y=\{y_1,\ldots,y_n\}$ and $\pa(v)=\{p_1,\ldots,p_n\}$, and define
the matrix $\AA$ as
\[
\AA_{ij}= \cases{ %
\bigl[(I-
\L)^T\Sigma\bigr]_{y_ip_j},&\quad $y_i\in\htr(v),$
\vspace*{2pt}\cr
\Sigma_{y_ip_j},&\quad $y_i\notin\htr(v).$}
\]
If $Y$ satisfies the half-trek criterion with respect to $v$, then
$\AA$ is generically invertible.
\end{lemma}
\begin{pf}
Recall the trek-rule from~(\ref{eqtrek-rule}). Let
$\mathcal{H}(v,w)\subset\mathcal{T}(v,w)$ be the set of all
half-treks from $v$ to $w$. Then, for each $i,j\in\{1,\ldots,n\}$,
\[
\AA_{ij}= \cases{ \displaystyle\sum_{\pi\in\mathcal{H}(y_i,p_j)} \pi(\lambda,
\omega), &\quad $y_i\in\htr(v),$\vspace*{2pt}
\cr
\displaystyle\sum
_{\pi\in\mathcal{T}(y_i,p_j)} \pi(\lambda,\omega), &\quad $y_i\notin
\htr(v).$}
\]
For a system of treks $\Pi$, define the monomial
\[
\Pi(\lambda,\omega) = \prod_{\pi\in\Pi} \pi(\lambda,
\omega).
\]
Then
\[
\det(\AA)=\sum_{\Psi\dvtx Y\rightrightarrows P} (-1)^{|\Psi|}\Psi(
\lambda,\omega),
\]
where the sum is over systems of treks $\Psi$ for which all treks
$\psi\in\Psi$ with sources in $ \htr(v)$ are half-treks. (The sign
$|\Psi|$ is the sign of the permutation that writes
$p_1,\ldots,p_n$ in the order of their appearance as targets of
the treks in $\Psi$.)

By assumption, there exists some system of half-treks with no sided
intersection from $Y$ to $P$. Let $\Pi$ be such a system, with
minimal total length among all such systems. Now take any system of
treks $\Psi$ from $Y$ to $P$, such that
$\Pi(\lambda,\omega)=\Psi(\lambda,\omega)$. (We do not assume that
$\Psi$ has no sided intersection, or has any half-treks.) In
Lemma~1 in the supplement\vadjust{\goodbreak} [\citet{Supplement}], we
prove that $\Psi=\Pi$
for any such $\Psi$---that is, $\Pi$ is the unique system of
half-treks with no sided intersection
of minimal total length. Therefore, the coefficient of the monomial
$\Pi(\lambda,\omega)$ in $\det(\AA)$ is given by $(-1)^{|\Pi|}$, and
$\det(\AA)$ is not the zero polynomial/power series. For generic
choices of $(\Lambda,\Omega)$ it thus holds that $\det(\AA)\not=0$.
\end{pf}

We now turn to the proof of the nonidentifiability theorem.

%
\begin{reptheorem}[(HTC-nonidentifiability)]
Suppose $G$ is a mixed graph in which every family $(Y_v\dvtx v\in
V)$ of subsets of the vertex set $V$ either contains a~set $Y_v$
that fails to satisfy the half-trek criterion with respect to $v$ or
contains a pair of sets $(Y_v,Y_w)$ with $v\in Y_w$ and
$w\in Y_v$. Then the parametrization $\phi_G$ is generically
infinite-to-one.
\end{reptheorem}

\begin{pf}
Let
\[
N= \bigl\{\{v,w\}\dvtx v\neq w,(v,w)\notin B \bigr\}
\]
be the set of (unordered) nonsibling pairs in the graph. Treating
$\Sigma$ as fixed, let $\JJ\in\mathbb{R}^{|N|\times|D|}$ be the
Jacobian of the equations in~(\ref{eqnzeros}), taking partial
derivatives with respect to the nonzero entries of $\Lambda$. The
entries of $\JJ$ are given by
%
\begin{equation}
\label{eqJac} \JJ_{\{v,w\},(u,v)}=- \bigl[(I-\L)^T\S
\bigr]_{wu} ,\qquad  \{v,w\} \in N, u\in\pa(v),
\end{equation}
and all other entries zero. By Lemma~2 in the
supplement, it is
sufficient to show that, under the conditions of the theorem, $\JJ$
does not have generically full column rank.

In the remainder of this proof, we always let
$\Sigma=\phi_G(\Lambda,\Omega)$ when considering~$\JJ$. If $\JJ$
has generically full column rank, then we can choose a set
$M\subset N$ with $|M|=|D|=\sum_{v\in V} |\pa(v)|$, such that
$\det(\JJ_{M,D})$ is not the zero polynomial, where $\JJ_{M,D}$ is the
square submatrix formed by taking all rows of $\JJ$ that are indexed
by $M$. By the definition of the determinant, there must be a
partition of $M=\bigcup_v M_v$ such that for all $v$, we have
\[
\det(\JJ_{M_v,(\pa(v),v)} )\neq0 .
\]
By~(\ref{eqJac}), each entry $\{w_1,w_2\}\in M_v$ must have
either $w_1=v$ or $w_2=v$. Writing $Y_v= \{w\dvtx\{v,w\}\in
M_v \}$, it holds that
\[
\det\bigl( \bigl[(I-\L)^T\S\bigr]_{Y_v,\pa(v)} \bigr)=\pm\det
(J_{\{Y_v,
v\},(\pa(v),v)} )= \pm\det(J_{M_v,(\pa(v),v)} )
\]
is nonzero. By Lemma~\ref{leminvert2} below, this implies that
each set $Y_v$ satisfies the half-trek criterion with respect to
its indexing node $v$. Forming a partition of $M\subset N$, the
sets $M_v$ are pairwise disjoint. Hence, no two nodes $v,w$ can
satisfy both $v\in Y_w$ and $w\in Y_v$ because otherwise
$\{v,w\}\in M_v\cap M_w$.
\end{pf}

%
\begin{lemma}
\label{leminvert2}
Let $v\in V$ be any node. Let $Y\subset V\setminus(\{v\}\cup
\sib(v) )$, with $|Y|=|\pa(v)|=n$. If the matrix
$\JJ=[(I-\L)^T\S]_{Y,\pa(v)}$ is generically invertible, then~$Y$
satisfies the half-trek criterion with respect to $v$.\vadjust{\goodbreak}
\end{lemma}
\begin{pf}
Abbreviate $P=\pa(v)$. We have
$\JJ=[(I-\L)^T\S]_{Y,P}=[\O(I-\L)^{-1}]_{Y,P}$. Hence,
\[
\det(\JJ)=\sum_{W\subset V,|W|=n} \det(\O_{Y,W})\det
\bigl((I-\L)^{-1}_{W,P}\bigr) .
\]
By assumption, $\det(\JJ)$ is not the zero polynomial/power series.
Therefore, for some $W\subset V$ with $|W|=n$, we have
$\det(\O_{Y,W})\not\equiv0$ and \mbox{$\det((I-\L)^{-1}_{W,P})\not
\equiv
0$}.

By Menger's theorem [see, e.g., Theorem
9.1 of \citet{schrijver2004a}], the nonvanishing of
$\det((I-\L)^{-1}_{W,P})$ implies that there is a system $\Psi$ of
pairwise vertex-disjoint directed paths
$\psi_i\dvtx w_i\dir\cdots\dir p_i$, $i\in[n]$, whose sources and targets
give $W=\{w_1,\ldots,w_n\}$ and $P=\{p_1,\ldots,p_n\}$, respectively.
Indeed, if no such system exists, then by Menger's theorem there is
a set $C$ of strictly less than $n$ vertices such that all directed paths
from $W$ to $P$ pass through~$C$. But this implies that the matrix
$(I-\L)^{-1}_{W,P}$ factors as $(I-\L)^{-1}_{W,C} \cdot
(I-\L)^{-1}_{C,P}$, and $|C|<n$ implies that
$\det((I-\L)^{-1}_{W,P})=0$, a contradiction. Note that
by erasing loops, we can further arrange that the $\psi_i$ do not have
self-intersections.

Since $\det(\O_{Y,W})\neq0$, we can index $Y=\{y_1,\ldots,y_n\}$
such that $\O_{y_iw_i}\neq0$ for all~$i$. This implies that either
$y_i=w_i$ or $y_i\bi w_i\in B$. Now define a system of half-treks
$\Pi\dvtx Y\rightrightarrows P$ by setting $\pi_i=\psi_i$ if $w_i=y_i$,
and extending $\psi_i$ at the left-hand side to
\[
\pi_i= y_i\bidir w_i\dir\cdots\dir
p_i
\]
if $y_i\neq w_i$. Since $\Psi$ has no sided intersection, $\Pi$
also has no sided intersection. It follows that $Y$ satisfies the
half-trek criterion with respect to $v$.
\end{pf}


\section{Conclusion}
\label{secconclusion}

We have proposed graphical criteria for determining identifiability as
well as nonidentifiability of linear structural equation models. The
criteria can be checked in time that is polynomial in the size of the
mixed graph representing the model. To our knowledge, they are the
best known. In particular, they apply to cyclic graphs. For acyclic
graphs, the graph decomposition method discussed in
Section~\ref{secgraph-decompositions} further extends their scope.
We expect the decomposition method to also extend the scope of the
criteria for cyclic graphs, when a cyclic model is suitably embedded
into an acyclic one, but we leave a thorough study of this problem for
future work.

Our algebraic computations revealed that there remains a ``gap''
between the necessary and the sufficient condition for rational
identifiability that we have developed. To better understand this
gap, it would be helpful to find an interesting class of graphs,
defined on an arbitrary number of nodes $m$, which is rationally
identifiable but not HTC-identifiable.

In models that are not HTC-identifiable, the half-trek method can
still prove certain parameters to be rationally identifiable; recall,
for instance, the example from Figure~\ref{figfiniteidexamples}(a).\vadjust{\goodbreak}
Referring to Theorem~\ref{mainthm}, if a set $Y_v$ satisfies the
half-trek criterion with respect to the indexing node $v$, and
$Y_v\cap\htr(v)=\varnothing$, then the proof of Theorem~\ref{mainthm}
shows how to obtain rational expressions in the covariance matrix
$\Sigma$ that equal the coefficients $\lambda_{wv}$, where $w\in
\pa(v)$. In the next step of the recursive procedure that proves
Theorem~\ref{mainthm}, we can solve for any node $u$ with $Y_u\cap
\htr(u)\subseteq\{v\}$. Continuing in this way, individual parameters
can be
identified even though ultimately the procedure will stop before all
nodes are visited, as we are discussing an HTC-inconclusive graph.
In particular, the maximum flow construction given in Algorithm \ref
{algHTC-id}
will reveal all nodes whose set of incoming directed edge parameters
can be identified via the half-trek criterion.
It would be interesting to compare this partial application of the
half-trek method to other graphical criteria for identification of
individual edge coefficients; see, in particular, \citet{garcia2010} for
a review and examples of such methods.

\section*{Acknowledgments}
This collaboration was started at a workshop at the American Institute
of Mathematics. We are grateful to Ilya Shpitser and Jin Tian for
helpful comments about existing literature.

\begin{supplement}
\stitle{Inconclusive graphs, proofs and algorithms}
\slink[doi]{10.1214/12-AOS1012SUPP} 
\sdatatype{.pdf}
\sfilename{aos1012\_supp.pdf}
\sdescription{The supplement starts with lists of some mixed graphs on
$m=5$ nodes that are not classifiable using our methods, to illustrate
the existing ``gap'' between our two criteria. After that we prove
lemmas used in the main paper for establishing the HTC-identifiability
and HTC-infinite-to-one criteria, and we provide details for the
results relating HTC-identifiability to GC-identifiability and to graph
decomposition. We then give correctness proofs for our algorithms for
checking the HTC-criteria, and we discuss the weak HTC-criteria. The
supplementary article concludes with a~computational-algebraic discussion
of the polynomial equations that led to the HTC-criteria.}
\end{supplement}

%


\printaddresses


\begin{thebibliography}{29}

\bibitem[\protect\citeauthoryear{Bollen}{1989}]{bollen1989}
\begin{bbook}[mr]
\bauthor{\bsnm{Bollen},~\bfnm{Kenneth~A.}\binits{K.~A.}}
(\byear{1989}).
\btitle{Structural Equations with Latent Variables}.
\bpublisher{Wiley}, \baddress{New York}.
\bid{mr={0996025}}
\bptok{imsref}%
\end{bbook}
\endbibitem

\bibitem[\protect\citeauthoryear{Brito}{2004}]{brito2004}
\begin{bmisc}[author]
\bauthor{\bsnm{Brito},~\bfnm{Carlos}\binits{C.}}
(\byear{2004}).
\bhowpublished{Graphical methods for identification in structural equation
  models. Ph.D. thesis, UCLA Computer Science Dept.}
\bptok{imsref}%
\end{bmisc}
\endbibitem

\bibitem[\protect\citeauthoryear{Brito and Pearl}{2002a}]{brito2002}
\begin{barticle}[mr]
\bauthor{\bsnm{Brito},~\bfnm{Carlos}\binits{C.}} \AND
  \bauthor{\bsnm{Pearl},~\bfnm{Judea}\binits{J.}}
(\byear{2002}a).
\btitle{A new identification condition for recursive models with correlated
  errors}.
\bjournal{Struct. Equ. Model.}
\bvolume{9}
\bpages{459--474}.
\bid{doi={10.1207/S15328007SEM0904_1}, issn={1070-5511}, mr={1930449}}
\bptok{imsref}%
\end{barticle}
\endbibitem

\bibitem[\protect\citeauthoryear{Brito and Pearl}{2002b}]{brito2002g}
\begin{binproceedings}[author]
\bauthor{\bsnm{Brito},~\bfnm{Carlos}\binits{C.}} \AND
  \bauthor{\bsnm{Pearl},~\bfnm{Judea}\binits{J.}}
(\byear{2002}b).
\btitle{A graphical criterion for the identification of causal effects in
  linear models}.
In \bbooktitle{Proceedings of the Eighteenth National Conference on Artificial
  Intelligence (AAAI)}
\bpages{533--538}.
\bpublisher{AAAI press}, \baddress{Palo Alto, CA}.
\bptok{imsref}%
\end{binproceedings}
\endbibitem

\bibitem[\protect\citeauthoryear{Brito and Pearl}{2006}]{brito2006}
\begin{binproceedings}[author]
\bauthor{\bsnm{Brito},~\bfnm{Carlos}\binits{C.}} \AND
  \bauthor{\bsnm{Pearl},~\bfnm{Judea}\binits{J.}}
(\byear{2006}).
\btitle{Graphical condition for identification in recursive {SEM}}.
In \bbooktitle{Proceedings of the Twenty-Second Conference on Uncertainty in
  Artificial Intelligence}
(\beditor{\bfnm{Rina}\binits{R.}~\bsnm{Dechter}} \AND
  \beditor{\bfnm{Thomas~S.}\binits{T.~S.}~\bsnm{Richardson}}, eds.)
\bpages{47--54}.
\bpublisher{AUAI Press}, \baddress{Arlington, VA}.
\bptok{imsref}%
\end{binproceedings}
\endbibitem

\bibitem[\protect\citeauthoryear{Chan and Kuroki}{2010}]{chan2010}
\begin{binproceedings}[author]
\bauthor{\bsnm{Chan},~\bfnm{Hei}\binits{H.}} \AND
  \bauthor{\bsnm{Kuroki},~\bfnm{Manabu}\binits{M.}}
(\byear{2010}).
\btitle{Using descendants as instrumental variables for the identification of
  direct causal effects in linear {SEM}s}.
In \bbooktitle{Proceedings of the Thirteenth International Conference on
  Artificial Intelligence and Statistics}
(\beditor{\bfnm{Yee~Whye}\binits{Y.~W.}~\bsnm{Teh}} \AND
  \beditor{\bfnm{Mike}\binits{M.}~\bsnm{Titterington}}, eds.).
\bseries{J. Mach. Learn. Res. (JMLR), Workshop and Conference Proceedings}
\bvolume{9}
\bpages{73--80}.
\bnote{Available at \url{http://jmlr.csail.mit.edu/proceedings/}.}
\bptok{imsref}%
\end{binproceedings}
\endbibitem

\bibitem[\protect\citeauthoryear{Cormen et~al.}{2001}]{cormen2001}
\begin{bbook}[mr]
\bauthor{\bsnm{Cormen},~\bfnm{Thomas~H.}\binits{T.~H.}},
  \bauthor{\bsnm{Leiserson},~\bfnm{Charles~E.}\binits{C.~E.}},
  \bauthor{\bsnm{Rivest},~\bfnm{Ronald~L.}\binits{R.~L.}} \AND
  \bauthor{\bsnm{Stein},~\bfnm{Clifford}\binits{C.}}
(\byear{2001}).
\btitle{Introduction to Algorithms}, \bedition{2nd} ed.
\bpublisher{MIT Press}, \baddress{Cambridge, MA}.
\bid{mr={1848805}}
\bptok{imsref}%
\end{bbook}
\endbibitem

\bibitem[\protect\citeauthoryear{Cox, Little and O'Shea}{2007}]{cox2007}
\begin{bbook}[mr]
\bauthor{\bsnm{Cox},~\bfnm{David}\binits{D.}},
  \bauthor{\bsnm{Little},~\bfnm{John}\binits{J.}} \AND
  \bauthor{\bsnm{O'Shea},~\bfnm{Donal}\binits{D.}}
(\byear{2007}).
\btitle{Ideals, Varieties, and Algorithms}, \bedition{3rd} ed.
\bpublisher{Springer}, \baddress{New York}.
\bid{doi={10.1007/978-0-387-35651-8}, mr={2290010}}
\bptok{imsref}%
\end{bbook}
\endbibitem

\bibitem[\protect\citeauthoryear{Decker et~al.}{2011}]{Singular}
\begin{bmisc}[author]
\bauthor{\bsnm{Decker},~\bfnm{Wolfram}\binits{W.}},
  \bauthor{\bsnm{Greuel},~\bfnm{Gert-Martin}\binits{G.-M.}},
  \bauthor{\bsnm{Pfister},~\bfnm{Gerhard}\binits{G.}} \AND
  \bauthor{\bsnm{Sch{\"o}nemann},~\bfnm{Hans}\binits{H.}}
(\byear{2011}).
\bhowpublished{\textsc{Singular} {3-1-3}---{A}~computer algebra system for
  polynomial computations. Available at
  \texttt{\href{http://www.singular.uni-kl.de}{http://}
  \href{http://www.singular.uni-kl.de}{www.singular.uni-kl.de}}.}
\bptok{imsref}%
\end{bmisc}
\endbibitem

\bibitem[\protect\citeauthoryear{Didelez, Meng and Sheehan}{2010}]{didelez2010}
\begin{barticle}[mr]
\bauthor{\bsnm{Didelez},~\bfnm{Vanessa}\binits{V.}},
  \bauthor{\bsnm{Meng},~\bfnm{Sha}\binits{S.}} \AND
  \bauthor{\bsnm{Sheehan},~\bfnm{Nuala~A.}\binits{N.~A.}}
(\byear{2010}).
\btitle{Assumptions of {IV} methods for observational epidemiology}.
\bjournal{Statist. Sci.}
\bvolume{25}
\bpages{22--40}.
\bid{doi={10.1214/09-STS316}, issn={0883-4237}, mr={2741813}}
\bptok{imsref}%
\end{barticle}
\endbibitem

\bibitem[\protect\citeauthoryear{Drton, Foygel and Sullivant}{2011}]{drton2011}
\begin{barticle}[mr]
\bauthor{\bsnm{Drton},~\bfnm{Mathias}\binits{M.}},
  \bauthor{\bsnm{Foygel},~\bfnm{Rina}\binits{R.}} \AND
  \bauthor{\bsnm{Sullivant},~\bfnm{Seth}\binits{S.}}
(\byear{2011}).
\btitle{Global identifiability of linear structural equation models}.
\bjournal{Ann. Statist.}
\bvolume{39}
\bpages{865--886}.
\bid{doi={10.1214/10-AOS859}, issn={0090-5364}, mr={2816341}}
\bptok{imsref}%
\end{barticle}
\endbibitem

\bibitem[\protect\citeauthoryear{Evans and Ringel}{1999}]{evans1999}
\begin{barticle}[author]
\bauthor{\bsnm{Evans},~\bfnm{William~N.}\binits{W.~N.}} \AND
  \bauthor{\bsnm{Ringel},~\bfnm{Jeanne~S.}\binits{J.~S.}}
(\byear{1999}).
\btitle{Can higher cigarette taxes improve birth outcomes?}
\bjournal{Journal of Public Economics}
\bvolume{72}
\bpages{135--154}.
\bptok{imsref}%
\end{barticle}
\endbibitem

\bibitem[\protect\citeauthoryear{Ford and Fulkerson}{1962}]{ford1962}
\begin{bbook}[mr]
\bauthor{\bsnm{Ford},~\bfnm{L.~R.}\binits{L.~R.} \bsuffix{Jr.}} \AND
  \bauthor{\bsnm{Fulkerson},~\bfnm{D.~R.}\binits{D.~R.}}
(\byear{1962}).
\btitle{Flows in Networks}.
\bpublisher{Princeton Univ. Press}, \baddress{Princeton, NJ}.
\bid{mr={0159700}}
\bptok{imsref}%
\end{bbook}
\endbibitem

\bibitem[\protect\citeauthoryear{Foygel, Draisma and Drton}{2012}]{Supplement}
\begin{bmisc}[author]
\bauthor{\bsnm{Foygel},~\bfnm{Rina}\binits{R.}},
  \bauthor{\bsnm{Draisma},~\bfnm{Jan}\binits{J.}} \AND
  \bauthor{\bsnm{Drton},~\bfnm{Mathias}\binits{M.}}
(\byear{2012}).
\bhowpublished{Supplement to ``Half-trek criterion for generic identifiability
  of linear structural equation models.''
  DOI:\doiurl{110.1214/12-AOS1012SUPP}.}
\bptok{imsref}%
\end{bmisc}
\endbibitem

\bibitem[\protect\citeauthoryear{Garcia-Puente, Spielvogel and
  Sullivant}{2010}]{garcia2010}
\begin{binproceedings}[author]
\bauthor{\bsnm{Garcia-Puente},~\bfnm{Luis~D.}\binits{L.~D.}},
  \bauthor{\bsnm{Spielvogel},~\bfnm{Sarah}\binits{S.}} \AND
  \bauthor{\bsnm{Sullivant},~\bfnm{Seth}\binits{S.}}
(\byear{2010}).
\btitle{Identifying causal effects with computer algebra}.
In \bbooktitle{Proceedings of the Twenty-sixth Conference on Uncertainty in Artificial
  Intelligence (UAI)}
(\beditor{\bfnm{Peter}\binits{P.}~\bsnm{Gr{\"{u}}nwald}} \AND
  \beditor{\bfnm{Peter}\binits{P.}~\bsnm{Spirtes}}, eds.).
\bpublisher{AUAI Press}. 
\bptok{imsref}%
\end{binproceedings}
\endbibitem

\bibitem[\protect\citeauthoryear{MathWorks Inc.}{2010}]{MATLAB2010}
\begin{bmisc}[author]
\borganization{MathWorks Inc.}
(\byear{2010}).
\bhowpublished{MATLAB version 7.10.0 (R2010a). Natick, MA.}
\bptok{imsref}%
\end{bmisc}
\endbibitem

\bibitem[\protect\citeauthoryear{Okamoto}{1973}]{okamoto1973}
\begin{barticle}[mr]
\bauthor{\bsnm{Okamoto},~\bfnm{Masashi}\binits{M.}}
(\byear{1973}).
\btitle{Distinctness of the eigenvalues of a quadratic form in a multivariate
  sample}.
\bjournal{Ann. Statist.}
\bvolume{1}
\bpages{763--765}.
\bid{issn={0090-5364}, mr={0331643}}
\bptok{imsref}%
\end{barticle}
\endbibitem

\bibitem[\protect\citeauthoryear{Pearl}{2000}]{pearl2000}
\begin{bbook}[mr]
\bauthor{\bsnm{Pearl},~\bfnm{Judea}\binits{J.}}
(\byear{2000}).
\btitle{Causality: Models, Reasoning, and Inference}.
\bpublisher{Cambridge Univ. Press}, \baddress{Cambridge}.
\bid{mr={1744773}}
\bptok{imsref}%
\end{bbook}
\endbibitem

\bibitem[\protect\citeauthoryear{R Development Core Team}{2011}]{R}
\begin{bmisc}[author]
\borganization{R Development Core Team}.
(\byear{2011}).
\bhowpublished{\textit{R: A language and environment for statistical computing}. R
  Foundation for Statistical Computing, Vienna, Austria.}
\bptok{imsref}%
\end{bmisc}
\endbibitem

\bibitem[\protect\citeauthoryear{Richardson and Spirtes}{2002}]{richardson2002}
\begin{barticle}[mr]
\bauthor{\bsnm{Richardson},~\bfnm{Thomas}\binits{T.}} \AND
  \bauthor{\bsnm{Spirtes},~\bfnm{Peter}\binits{P.}}
(\byear{2002}).
\btitle{Ancestral graph {M}arkov models}.
\bjournal{Ann. Statist.}
\bvolume{30}
\bpages{962--1030}.
\bid{doi={10.1214/aos/1031689015}, issn={0090-5364}, mr={1926166}}
\bptok{imsref}%
\end{barticle}
\endbibitem

\bibitem[\protect\citeauthoryear{Schrijver}{2004}]{schrijver2004a}
\begin{bbook}[author]
\bauthor{\bsnm{Schrijver},~\bfnm{Alexander}\binits{A.}}
(\byear{2004}).
\btitle{Combinatorial Optimization. {P}olyhedra and Efficiency}.
\bseries{Algorithms and Combinatorics 24}
\bvolume{A}.
\bpublisher{Springer}, \baddress{Berlin}.
\bptok{imsref}%
\end{bbook}
\endbibitem

\bibitem[\protect\citeauthoryear{Spirtes, Glymour and
  Scheines}{2000}]{spirtes2000}
\begin{bbook}[mr]
\bauthor{\bsnm{Spirtes},~\bfnm{Peter}\binits{P.}},
  \bauthor{\bsnm{Glymour},~\bfnm{Clark}\binits{C.}} \AND
  \bauthor{\bsnm{Scheines},~\bfnm{Richard}\binits{R.}}
(\byear{2000}).
\btitle{Causation, Prediction, and Search}, \bedition{2nd} ed.
\bpublisher{MIT Press}, \baddress{Cambridge, MA}.
\bid{mr={1815675}}
\bptok{imsref}%
\end{bbook}
\endbibitem

\bibitem[\protect\citeauthoryear{Sullivant, Talaska and
  Draisma}{2010}]{sullivant2010}
\begin{barticle}[mr]
\bauthor{\bsnm{Sullivant},~\bfnm{Seth}\binits{S.}},
  \bauthor{\bsnm{Talaska},~\bfnm{Kelli}\binits{K.}} \AND
  \bauthor{\bsnm{Draisma},~\bfnm{Jan}\binits{J.}}
(\byear{2010}).
\btitle{Trek separation for {G}aussian graphical models}.
\bjournal{Ann. Statist.}
\bvolume{38}
\bpages{1665--1685}.
\bid{doi={10.1214/09-AOS760}, issn={0090-5364}, mr={2662356}}
\bptok{imsref}%
\end{barticle}
\endbibitem

\bibitem[\protect\citeauthoryear{Tian}{2005}]{tian2005}
\begin{binproceedings}[author]
\bauthor{\bsnm{Tian},~\bfnm{Jin}\binits{J.}}
(\byear{2005}).
\btitle{Identifying direct causal effects in linear models}.
In \bbooktitle{Proceedings of the Twentieth National Conference on Artificial
  Intelligence (AAAI)}
\bpages{346--353}.
\bpublisher{AAAI press}, \baddress{Palo Alto, CA}.
\bptok{imsref}%
\end{binproceedings}
\endbibitem

\bibitem[\protect\citeauthoryear{Tian}{2009}]{tian2009}
\begin{binproceedings}[author]
\bauthor{\bsnm{Tian},~\bfnm{Jin}\binits{J.}}
(\byear{2009}).
\btitle{Parameter identification in a class of linear structural equation
  models}.
In \bbooktitle{Proceedings of the Twenty-first International Joint Conference on Artificial
  Intelligence (IJCAI)}
\bpages{1970--1975}.
\bpublisher{AAAI press}, \baddress{Palo Alto, CA}.
\bptok{imsref}%
\end{binproceedings}
\endbibitem

\bibitem[\protect\citeauthoryear{Wermuth}{2011}]{wermuth2010}
\begin{barticle}[mr]
\bauthor{\bsnm{Wermuth},~\bfnm{Nanny}\binits{N.}}
(\byear{2011}).
\btitle{Probability distributions with summary graph structure}.
\bjournal{Bernoulli}
\bvolume{17}
\bpages{845--879}.
\bid{doi={10.3150/10-BEJ309}, issn={1350-7265}, mr={2817608}}
\bptok{imsref}%
\end{barticle}
\endbibitem

\bibitem[\protect\citeauthoryear{Wright}{1921}]{wright1921}
\begin{barticle}[author]
\bauthor{\bsnm{Wright},~\bfnm{Sewall}\binits{S.}}
(\byear{1921}).
\btitle{Correlation and causation}.
\bjournal{J. Agricultural Research}
\bvolume{20}
\bpages{557--585}.
\bptok{imsref}%
\end{barticle}
\endbibitem

\bibitem[\protect\citeauthoryear{Wright}{1934}]{wright1934}
\begin{barticle}[author]
\bauthor{\bsnm{Wright},~\bfnm{Sewall}\binits{S.}}
(\byear{1934}).
\btitle{The method of path coefficients}.
\bjournal{Ann. Math. Statist.}
\bvolume{5}
\bpages{161--215}.
\bptok{imsref}%
\end{barticle}
\endbibitem

\end{thebibliography}
\end{document}